\documentclass[12pt,a4paper]{amsart}
\usepackage{geometry}
\usepackage{amsmath, amsthm, amssymb, amsfonts}
\usepackage[pdftex,bookmarks=true]{hyperref}
\usepackage{color}
\usepackage[mathscr]{euscript}
\usepackage{enumitem}
\usepackage{todonotes}

\newtheorem{theorem}{Theorem}[section]
\newtheorem{proposition}{Proposition}[section]
\newtheorem{lemma}[theorem]{Lemma}

\newtheorem{observation}{Observation}[section]

\newtheorem{remark}[theorem]{Remark}

\def\C{{\mathbb C}}  % set of complex numbers
\def\Z{{\mathbb Z}} % set of integers

\def\R{{\mathbb R}}

\def\r{{\mathcal{R}}}
\def\P{{\mathbb P}}

\def\I{{\mathcal{I}}}

\def\supp{{\rm supp}}

\def\Z{{\mathbb Z}}

\def\<{\left<}
\def\>{\right>}

\def\dist{{\rm dist}}
\def\size{{\rm size\,}}
\def\Rest{{\rm Rest}}
\def\ave{{\rm \ ave}}

\newtheorem*{notation}{Notation:}

\newcommand{\one}{\mathbf{1}}
\newcommand{\rr}{\mathbb}
\newcommand{\ic}{\mathcal}
\newcommand{\ii}{\mathcal}
\newcommand{\ci}{\tilde{\chi}}

\newcommand{\ds}{\displaystyle}

\def\Xint#1{\mathchoice
   {\XXint\displaystyle\textstyle{#1}}%
   {\XXint\textstyle\scriptstyle{#1}}%
   {\XXint\scriptstyle\scriptscriptstyle{#1}}%
   {\XXint\scriptscriptstyle\scriptscriptstyle{#1}}%
   \!\int}
\def\XXint#1#2#3{{\setbox0=\hbox{$#1{#2#3}{\int}$}
     \vcenter{\hbox{$#2#3$}}\kern-.5\wd0}}

\def\aver#1{\Xint-_{#1}}

\title[Mixed norm estimates for Littlewood-Paley square functions]{Multiple vector-valued, mixed norm estimates for Littlewood-Paley square functions}

\author
[C. Benea, C. Muscalu]
{Cristina Benea\ \ \ Camil Muscalu}

\address{Universit\'{e} de Nantes,
Laboratoire Jean Leray, Nantes, 44311, France}
\email{cristina.benea@univ-nantes.fr}

\address{Department of Mathematics,
Cornell University, Ithaca, NY 14853, USA \footnote{The author is also a Member of the Simion Stoilow Institute of Mathematics of the Romanian Academy in Bucharest.}}
\email{camil@math.cornell.edu}


\begin{document}

\begin{abstract}
We prove that for any $L^Q$-valued Schwartz function $f$ defined on $\R^d$, one has the multiple vector-valued, mixed norm estimate
$$
\| f \|_{L^P(L^Q)} \lesssim \| S f \|_{L^P(L^Q)}
$$
valid for every $d$-tuple $P$ and every $n$-tuple $Q$ satisfying $0 < P, Q < \infty$ componentwise.  Here $S:= S_{d_1}\otimes ... \otimes S_{d_N}$ is a tensor product of several Littlewood-Paley
square functions $S_{d_j}$ defined on arbitrary Euclidean spaces $\R^{d_j}$ for $1\leq j\leq N$, with the property that $d_1 + ... + d_N = d$.
This answers a question that came up
implicitly in our recent works \cite{BM1}, \cite{BM2}, \cite{BM3} and completes in a natural way classical results of the Littlewood-Paley theory.
The proof is based on the {\it helicoidal method} introduced by the authors in the aforementioned papers.
\end{abstract}
\maketitle

\section{Introduction}
\label{sec:introduction}
Let us start by recalling that a  sequence of $L^1$-bounded Schwartz functions $(\psi_k)_{k\in\Z}$ defined on the Euclidean space $\R^m$  is called a Littlewood-Paley sequence, if its Fourier transform satisfies \footnote{Here and throughout the article we use the standard notation $A\lesssim B$ meaning that $A \leq C B$ for some constant $C > 0$
which can be universal, or depending on several implicit parameters, coming from the specific context.}
\begin{equation}
0\notin \supp \,\widehat{\psi_k} , \qquad |\partial^{\alpha} \widehat{\psi_k}(\xi)|  \lesssim 2^{-\alpha k} \big(1+ \frac{|\xi|}{2^k}\big)^{-100\, m}
\end{equation}
for every $\xi\in\R^m$ and sufficiently many multi-indices $\alpha$, and if one also has
$$1 = \sum_{k\in\Z} \widehat{\psi_k} .$$
In particular, any Schwartz function $f$ defined on $\R^m$ admits the following Littlewood-Paley decomposition
$$f = \sum_{k\in\Z} f\ast \psi_k .$$
To any Littlewood-Paley sequence, one can also associate a Littlewood-Paley square function $S_m f$, defined by
\begin{equation}
S_m f (x):= \big( \sum_{k\in \Z} |f \ast \psi_k (x) |^2 \big)^{1/2}.
\end{equation}
Moreover, for any $N\geq 1$ such Littlewood-Paley sequences $(\psi^j_k)_{k\in\Z}$ defined on $\R^{d_j}$ for $1\leq j\leq N$, one defines an $N$-parameter one $(\Psi_k)_{k\in\Z^N}$ on $\R^d:= \R^{d_1} \times ... \times \R^{d_N}$
by
\begin{equation}
\Psi_k := \psi^1_{k_1}\otimes ... \otimes \psi^N_{k_N}
\end{equation}
for $k = (k_1, ..., k_N)$, where
$$
\psi^1_{k_1}\otimes ... \otimes \psi^N_{k_N} (x_1, ..., x_N) := \psi^1_{k_1}(x_1)\cdot ... \cdot \psi^N_{k_N}(x_N).
$$
Here we think of the generic variable $x\in\R^d$ as being identified with the vector $(x_1, ..., x_N)$ with $x_j \in \R^{d_j}$ for $1 \leq j \leq N$.

In particular, any Schwartz function on $\R^d$ admits the decomposition

$$f = \sum_{k\in\Z^N} f\ast \Psi_k .$$
One can then also define the $N$-parameter square function $S f$ by the formula
\begin{equation}\label{square}
S f (x) : = \big( \sum_{k\in\Z^N} |f\ast \Psi_k (x) |^2 \big)^{1/2}
\end{equation}
for $x\in \R^d$. This is the square function that will be studied in the present article.

To complete the presentation of the main notations that we will use, we also recall that given any $n\geq 1$ {\color{black}$\sigma$}-finite measurable spaces $(S_j, \Sigma_j, \mu_j)$ for $1\leq j\leq n$ and $R=(r_1, ... , r_n)$ an $n$-tuple of positive real numbers, 
one can define the iterated (or mixed norm) Lebesgue
space $L^R(S, \Sigma, \mu)$ to be the space containing those functions $g$ which are measurable on the product space
$$(S, \Sigma, \mu) := (\prod_{j=1}^n S_j, \prod_{j=1}^n \Sigma_j, \prod_{j=1}^n \mu_j)$$
and for which the (quasi)-norm $\|g\|_R$ defined by
$$\|g\|_R : = \| ... \| g(s_1, ..., s_n)\|_{L^{r_n}(S_n, \Sigma_n, \mu_n)} ... \|_{L^{r_1}(S_1, \Sigma_1, \mu_1)}$$
is finite.

The classical Littlewood-Paley theory states that the following inequalities 
\begin{equation}\label{prima}
\|f\|_{L^p(\R^m)} \lesssim \| S_m f \|_{L^p(\R^m)} \lesssim \|f\|_{L^p(\R^m)}
\end{equation}
are true, provided that $1< p< \infty$ and that, in addition, the left hand side of (\ref{prima})
\begin{equation}\label{doua}
\|f\|_{L^p(\R^m)} \lesssim \| S_m f \|_{L^p(\R^m)} 
\end{equation}
is in fact available in the whole range  $0 < p < \infty$, see for instance \cite{cw} and \cite{stein}.

Standard duality and vector-valued arguments for singular integrals allow one to extend (\ref{prima}) very easily to the setting of mixed norm spaces and $N$-parameter square functions. This implies that the inequalities
 \begin{equation}\label{treia}
\|f\|_{L^P(L^Q)} \lesssim \| S f \|_{L^P(L^Q)} \lesssim \|f\|_{L^P(L^Q)}
\end{equation}
are true for $L^Q$-valued Schwartz functions defined in $\R^d$ for every $n$-tuple $Q$ and $d$-tuple $P$ satisfying $1 < P, Q < \infty$ componentwise.

To be more specific, the space $L^P$ above is considered with respect to the product Lebesgue measure in $\R^d$, and as before, by
$\|h\|_{L^P(L^Q)}$ one means the mixed (quasi)-norm given by
$$\|h\|_{L^P(L^Q)} : = \| \,\,\| h(x,s) \|_{L^Q (S,\, \Sigma, \, \mu)}\,\|_{L^P (\R^d)}.$$

The main result of the present article is to show that a similar extension can be proved for the estimate (\ref{doua}).

\begin{theorem}\label{mainth}
The following estimate  
\begin{equation}\label{mainest}
\|f\|_{L^P(L^Q)} \lesssim \| S f \|_{L^P(L^Q)}
\end{equation}
is true, for every $L^Q$-valued Schwartz function $f$ on $\R^d$, as long as the $n$-tuples $Q$ and the $d$-tuples $P$ satisfy the condition $0 < P, Q < \infty$ componentwise.

\end{theorem}

As we will see, unlike (\ref{treia}), the proof of Theorem \ref{mainth} is far from being routine, and it is based on the {\it helicoidal method} developed by the authors in \cite{BM1}, \cite{BM2}, \cite{BM3} .    
The question addressed and answered by Theorem \ref{mainth} surfaced out quite naturally in our recent works  \cite{BM1}, \cite{BM2} and it is related to an open problem of Kenig on mixed norm estimates
for paraproducts on polydisks. See also our recent expository work \cite{BM4}, in particular Theorem 5 there.

Some particular cases of (\ref{mainest}) were known in the scalar case, that is when $L^Q = \C$.  The case when all the entries of the $d$-tuple $P$ are equal to each other is the well known multi-parameter case studied by Gundy and Stein in \cite{GS}.
More recently, Hart, Torres and Wu have proved the case when  $N=1$ and $d=2$, again, in the scalar situation \cite{hrw}. Even more recently, the case $N=1$ was extended in \cite{anisHardySpMixedNormLP} to arbitrary dimensions, also in an anisotropic setting.

The central point of the paper will be the proof of our main Theorem \ref{mainth} based on techniques from \cite{BM1}, \cite{BM2} and \cite{BM3}. We split the presentation into two distinct parts. In the first part, we consider the case when all the square functions $S_{d_j}$ for $1\leq j\leq N$ are one
dimensional, that is when $d_1= ... = d_N =1$. Notice that in this case $N=d$. The proof of this case represents the core of the present article.

Under this assumption, we first show in Section \ref{sec:reduction}, that the estimate (\ref{mainest}) follows easily, by induction, from its particular case $d=1$. Notice that in this situation, (\ref{mainest}) becomes a multiple vector-valued extension of the well known (scalar) inequality (\ref{doua}). Then, in Section \ref{sec:discrete-mvv-case}, we explain how this multiple vector-valued case is implied by a certain discrete analogue of it.

 Next, in Section \ref{sec:proof-discretized-main-thm}, which is more involved, we describe the proof of this discrete case, by using ideas that lie at the heart of our helicoidal method in \cite{BM1}, \cite{BM2}, \cite{BM3}. In Section \ref{sec:proof-main-thm-gen-case} we explain how one can modify the proof in part one to handle the general,
mixed norm case, of Theorem \ref{mainth}.

Lastly, in the final Section \ref{sec:weights-extrap}, we will see how Theorem \ref{mainth} can also be obtained through extrapolation from a weighted, scalar version of Theorem \ref{mainth}, which appeared in the context of weighted Hardy spaces in \cite{DingHanLu-weightedHardySpaces}. Since we are outside the Banach setting, the extrapolation needed concerns $A_\infty$ weights and pairs of functions. For the mixed-norm estimates, we need to adapt a result of Kurtz \cite{KurtzExtrapMixedNormedSpaces}.

That the vector-valued result of Theorem \ref{mainth} allows also for a proof based on extrapolation and weighted theory should not be surprising: the helicoidal method yields vector-valued results that can be obtained also through extrapolation, once weighted estimates for the \emph{correct} class of weights is known. This was the case also with the bilinear Hilbert transform (see \cite{BM1}, \cite{BM3}, \cite{extrap-BHT}, \cite{martell-kangwei-mulilinear-weights-extrapolation}). For completeness, in Section \ref{sec:weighted-result-hel-method} we show how to deduce the weighted version of Theorem \ref{mainth} by using the helicoidal method: the same maximal inequality used in Section \ref{sec:proof-discretized-main-thm} plays a central role, and only the stopping time algorithm changes.

{\bf Acknowledgements :} C.B. was partially supported by the ERC Project FAnFArE no. 637510. C.M. was partially supported by the NSF Grant DMS 1500262. He also acknowledges partial support from a grant from the Ministry of Research and Innovation of Romania, CNCS - UEFISCDI, project PN-III-P4-ID-PCE-2016-0823 within PNCDI - III. During the Spring Semester of 2017, C.M. was a member of the MSRI in Berkeley, as part of the Program in Harmonic Analysis, and during the Fall Semester of 2017, he was visiting the Mathematics Department of the Universit\'{e} Paris-Sud Orsay, as a Simons Fellow. He is grateful to both institutions for their hospitality, and to the Simons Foundation for their generous support.

The authors are grateful to Dachun Yang for pointing out the results in \cite{anisHardySpMixedNormLP}.

\section{Reduction to the multiple vector-valued case}
\label{sec:reduction}
As mentioned above, we first study the case when $d_1= ... = d_N =1$. From now on, until the last section of the paper, we work under this assumption.

And as also mentioned in the introduction, in this section we show that Theorem \ref{mainth} follows by induction, from its particular case $d=1$. Recall also that $d=N$ now. Let us therefore assume that Theorem \ref{mainth}
is true for dimensions smaller or equal than $d-1$ and we will explain how to deduce the $d$ dimensional case. The argument is based on the following identity

\begin{equation}\label{squared-1}
 S f = \big( \sum_{k \in \Z} | S_{(x_1, ..., x_{d-1})} (f\ast \psi^d_k)|^2 \big)^{1/2}
\end{equation}
where $S_{(x_1, ..., x_{d-1}) }(g)$ denotes the $d-1$ dimensional part of the square function, taken with respect to the variables $x_1, ..., x_{d-1}$, and explicitly given by

\begin{equation}\label{adoua}
S_{(x_1, ..., x_{d-1})} (g) (x) : = \Big( \sum_{k_1, ..., k_{d-1}} | g\ast (\psi^1_{k_1}\otimes  ... \otimes \psi^{d-1}_{k_{d-1}})(x) |^2 \Big)^{1/2}
\end{equation}
for $x\in\R^d$.  The first  convolution in (\ref{squared-1}) is a one dimensional one, taken with respect to the last variable $x_d$, while  the convolution in (\ref{adoua}) is a $d-1$ dimensional one, taken 
with respect to the first $d-1$ variables $x_1, ..., x_{d-1}$ . Using (\ref{squared-1}) one can write
\begin{align*}
\| S f \|_{L^P(L^Q)} & =  \big\| \big( \sum_{k\in\Z} | S_{(x_1, ..., x_{d-1})} (f\ast \psi^d_k)|^2 \big)^{1/2} \big\|_{L^P(L^Q)} \\
&= \big\| \big( S_{(x_1, ..., x_{d-1})}(f\ast \psi^d_k)\big)_k \big\|_{L^P(L^Q(\ell^2))} =  \big\| \big( S_{(x_1, ..., x_{d-1})}(f\ast \psi^d_k)\big)_k \big\|_{L^{\widetilde{P}}(L^{p_d}(L^Q(\ell^2)))}, 
\end{align*}
where $\widetilde{P} := (p_1, ..., p_{d-1})$.

Here, one can use the induction hypothesis in the $(d-1)$ dimensional case to conclude that the above expression is larger than

$$\big\| \big( f\ast \psi^d_k\big)_k \big\|_{L^{\widetilde{P}}(L^{p_d}(L^Q(l^2)))} = \big\| \big( \sum_k |f\ast\psi^d_k|^2\big)^{1/2}\big\|_{L^{\widetilde{P}}(L^{p_d}(L^Q))} .$$
Finally, by using the one dimensional case and Fubini, we see that this is also greater than
$$\|f \|_{L^{\widetilde{P}}(L^{p_d}(L^Q))} = \|f \|_{L^P(L^Q)},$$
which ends the argument.

\section{The discrete multiple vector-valued case}
\label{sec:discrete-mvv-case}
Now that we know that Theorem \ref{mainth} (in the special situation when $d_1= ... = d_N =1$) can be reduced to its $d=1$ particular case, we show in this section that a further reduction is possible. The multiple vector-valued $d=1$ case can be reduced to
a {\it discrete} variant of it, that will be described next.

Let us pause briefly and recall that a sequence of Schwartz functions $(\phi_I)_I$ on the real line, indexed by dyadic intervals $I$, is called an $L^p$ normalized {\it lacunary} sequence (for some $p\in (0, \infty]$), if and only if the following estimates hold

\begin{equation}\label{eq0}
|\partial^{\alpha} \phi_I(x)| \lesssim \frac{1}{|I|^{1/p}} \frac{1}{|I|^{\alpha}}\big(1+ \frac{\dist(x, I)}{|I|}\big)^{-100}
\end{equation}
for $x\in\R$, $0\leq \alpha \leq 10$ and also if $\int_{\R} \phi_I(x) d x = 0$.

Let now $(\phi^1_I)_I$ and $(\phi^2_I)_I$ be two $L^2$-normalized such lacunary sequences, indexed by a finite arbitrary subset of dyadic intervals. The following discrete variant of the one dimensional case of
Theorem \ref{mainth} is true.

\begin{theorem}\label{discreteth}
For every $0 < p < \infty$ and tuple $Q$ as before, one has
\begin{equation}\label{discreteest}
\| \sum_I \langle f, \phi^1_I\rangle \phi^2_I \|_{L^p(L^Q)} \lesssim \Big\| \big( \sum_I \frac{ |\langle f, \phi^1_I\rangle|^2  } {|I|} \one_I \big)^{1/2}\Big\|_{L^p(L^Q)}.
\end{equation}
\end{theorem}

\begin{observation}\label{obs}
The function $f$ above depends on the variables $(s_1, ..., s_n) \in S$ and on $x\in \R$. Sometimes we will write this explicitly as $f_{(s_1, ..., s_n)}(x)$. It is important to emphasize
that, as we will see from the proof of Theorem \ref{discreteth}, the estimate (\ref{discreteest} ) holds also in the more general case when the families $(\phi^1_I)_I$ and $(\phi^2_I)_I$ depend on the variables
$(s_1, ..., s_n)\in S$ as well, in a uniform manner, with respect to the implicit constants of (\ref{eq0}).
\end{observation}
We explain now why Theorem \ref{discreteth} implies the one dimensional case of Theorem \ref{mainth}. The argument is based on an idea that we learned from the article \cite{hrw}, and which goes back to the work of Frazier and Jawerth \cite{FJ}.  

\begin{proposition}\label{descpr}
There exists a large universal constant $N$ such that, given any sequence of intermediate points $x_I\in I$, there exists $(\psi_I)_I$  an $L^{\infty}$ normalized lacunary sequence, so that every Schwartz function $h$ on the real line can be decomposed as
\begin{equation}\label{desc}
h= \sum_k \sum_{|I| = 2^{-k}} (h\ast \psi_{k-N})(x_I) \psi_I.
\end{equation}
\end{proposition}
In (\ref{desc}), the sequence $(\psi_l)_l$ is any a priori fixed Littlewood-Paley sequence. We prove Proposition \ref{descpr} in detail later on. In what follows, we describe how it helps reducing the $d=1$ case of Theorem \ref{mainth} to its discrete analogue from Theorem \ref{discreteth}.

Fix $f (= f_{(s_1, ..., s_n)}(x) )$. For every $(s_1, ..., s_n) \in S$ pick $x_I\in I$ a number with the property that
$$
\inf_{y\in I} \left| (f_{(s_1, ..., s_n)}\ast \psi_{k-N})(y)\right| = \left| (f_{(s_1, ..., s_n)}\ast \psi_{k-N})(x_I)\right| 
$$
where $I$ is a dyadic interval with $|I| = 2^{-k}$. Clearly, $x_I$ depends on $f$ and also, implicitly, on $(s_1, ..., s_n)\in S$.

Using Proposition \ref{descpr}, one can write
\begin{equation}\label{est2}
\|f\|_{L^p(L^Q)} = \| \sum_k \sum_{|I| = 2^{-k}} (f_{(s_1, ..., s_n)}\ast \psi_{k-N})(x_I) \psi_{I, (s_1, ..., s_n)} (x) \|_{L^p(L^Q)}.
\end{equation}
Using now the general form of Theorem \ref{discreteth} (see Observation \ref{obs} that followed it) one can majorize the above expression (\ref{est2}) further by
$$
\Big\|\big(\sum_k \sum_{|I|=2^{-k}} | (f_{(s_1, ..., s_n)} \ast \psi_{k-N})(x_I) |^2 \one_I(x) \big)^{1/2}\Big\|_{L^p(L^Q)}
$$
and using the definition of the sequence $(x_I)_I$ above, one can immediately see that this is smaller than
$$
\Big\|\big(\sum_k  | (f_{(s_1, ..., s_n)} \ast \psi_{k-N})(x) |^2 \big)^{1/2}\Big\|_{L^p(L^Q)} =
\Big\| \big( \sum_k |f\ast \psi_k|^2\big)^{1/2}\Big\|_{L^p(L^Q)}$$
as desired.

\subsection{Proof of Proposition \ref{descpr}}

We now describe the proof of Proposition \ref{descpr} using the ideas from \cite{FJ}. 

Start by writing, for a generic function of one variable $f$ :

$$f = \sum_k f\ast\psi_k = \sum_k f\ast \psi_{k-N}.$$
We will prove that for every $k\in \Z$, a family of functions $(\psi_I)_I$ as in Proposition \ref{descpr} exists\footnote{This time all the intervals $I$ have the same length, $|I| = 2^{-k}$.},
so that
\begin{equation}\label{eq1}
f\ast \psi_{k-N} = \sum_{|I|=2^{-k}} (f\ast \psi_{k-N})(x_I) \psi_I .
\end{equation}
Clearly, this would be enough. Since the argument is scale invariant, we will prove this in the particular case when $k=N$. In this case, (\ref{eq1}) becomes

\begin{equation}\label{eq2}
f\ast \psi_0 = \sum_{|I| = 2^{-N}} (f\ast \psi_0)(x_I)  \psi_I .
\end{equation}
Consider now $\widetilde{\psi}$ a Schwartz function so that $\widehat{\widetilde{\psi}} = 1$ on the support of $\widehat{\psi_0}$ and having the property that
$\supp \, \widehat{\widetilde{\psi}} \subseteq [1/2, 4]$. 

Then, one can write
{\fontsize{10}{10}\begin{align*}
&f\ast \psi_0(x)= (f\ast \psi_0)\ast \widetilde{\psi}(x) = \int_{\R} f\ast \psi_0(y) \widetilde{\psi}(x-y) d y  = \sum_{|I| = 2^{-N}} \int_{\R} f\ast \psi_0(y) \widetilde{\psi}(x-y) \one_I(y) d y \\
&= \sum_{|I|=2^{-N}} f\ast\psi_0(x_I) \int_{\R} \widetilde{\psi}(x-y) \one_I(y) d y + \Rest_1(x)= \sum_{|I|=2^{-N}} f\ast\psi_0(x_I) \int_{\R} \widetilde{\psi}(x-y) \one_I(y) d y + \Rest_1(x),
\end{align*}}
where $\phi^1_I := \widetilde{\psi}\ast \one_I(x)$ and
\begin{equation}\label{eq4}
\Rest_1(x) = : \sum_{|I|=2^{-N}} \Rest_{1, I}(x)= \sum_{|I|=2^{-N}} \int_{\R} [ f\ast \psi_0(y) - f\ast\psi_0(x_I) ] \widetilde{\psi}(x-y) \one_I(y) d y .
\end{equation}

The above inner expression can be estimated by
\begin{align*}
f\ast\psi_0(y) - f \ast \psi_0(x_I) &= \int_{\R} f(z)  [ \psi_0(y-z) - \psi_0(x_I-z) ] d z  \\
& =  \int_{\R} f(z) \psi'_0 (\# - z) (y - x_I) d z,
\end{align*}
where $\#$ is a point lying inside the interval $I$ and depending on $y, x_I$ and $z$. Since both $y$ and $x_I$ belong to $I$, it is easy to see that the above expression
is at most $C \, 2^{-N} \| f\|_{\infty}$. Using this in (\ref{eq4}) we obtain that
$$|\Rest_{1,I}(x)| \leq C_M 2^{-N} \|f\|_{\infty}  (1+ \dist (x, I))^{-M} |I|$$
which implies further
$$|\Rest_1(x)| \leq C \|f\|_{\infty} 2^{-N}.$$
We see these calculations as providing a first approximation towards the desired (\ref{eq2}). To summarize, so far we showed that
\begin{equation}\label{eq5}
f\ast \psi_0(x) = \sum_{|I|=2^{-N}} f\ast\psi_0(x_I) \phi^1_I(x) + \Rest_1(x)
\end{equation}
where $|\Rest_1(x)| \leq C\,\|f\|_{\infty}2^{-N}$ and $(\phi^1_I)_I$ is a lacunary family.

We now iterate this fact carefully. Fix $J$ with $|J|=2^{-N}$ and recall the following expression
\begin{equation}\label{eq6}
\Rest_{1,J}(x) = \int_{\R} [ f\ast\psi_0(y) - f\ast\psi_0(x_J)] \widetilde{\psi}(x-y) \one_J(y) d y.
\end{equation}
Using (\ref{eq5}) for $x=y$ and $x=x_J$ in (\ref{eq6}) we obtain a decomposition of $\Rest_{1,J}(x)$ of type
\begin{align*}
&\sum_{|I|=2^{-N}} f\ast\psi_0(x_I) \int_{\R} [ \phi^1_I(y) - \phi^1_I(x_J)] \widetilde{\psi}(x-y) \one_J(y) d y \\
&+ \sum_{|I|=2^{-N}} \int_{\R} [ \Rest_{1,I}(y) - \Rest_{1,I}(x_J)] \widetilde{\psi}(x-y) \one_J(y) d y .
\end{align*}

Summing over $|J|=2^{-N}$, we obtain the formula
$$\Rest_1(x) = \sum_{|I|=2^{-N}} f\ast\psi_0(x_I) \phi^2_I + \Rest_2(x)$$
where
$$\phi^2_I(x) := \sum_{|J|=2^{-N}} \int_{\R} [ \phi^1_I(y) - \phi^1_I(x_J) ] \widetilde{\psi} (x-y) \one_J(y) d y$$
while
$$\Rest_2(x) = \sum_{|I|=2^{-N}} \Rest_{2,I}(x) $$
and
$$\Rest_{2,I}(x) := \sum_{|J|=2^{-N}} \int_{\R} [ \Rest_{1,I}(y) - \Rest_{1,I}(x_J) ] \widetilde{\psi}(x-y) \one_J(y) d y. $$
Arguing exactly as before, given that both $y$ and $x_J$ belong to the interval $J$, it is not difficult to see that $(\phi^2_I)_I$ is a lacunary family satisfying

$$\|\phi^2_I\|_{\infty} \leq C 2^{-N}, \quad \text{while}\quad \|\Rest_2\|_{\infty} \leq C^2 2^{-2N} \|f\|_{\infty}$$
where as always, $C$ is a universal constant. In other words, at our second approximation step, we obtain the decomposition
$$f \ast\psi_0(x) = \sum_{|I| = 2^{-N}} f\ast\psi_0(x_I) (\phi^1_I(x) + \phi^2_I(x)) + \Rest_2(x).$$
Iterating this an arbitrary number of times, we obtain that $f\ast\psi_0(x)$ can be written as

\begin{equation}\label{eq7}
f\ast\psi_0(x) = \sum_{|I|=2^{-N}} f\ast\psi_0(x_I) (\phi^1_I(x) + ... + \phi^l_I(x)) + \Rest_l(x)
\end{equation}
where $(\phi^j_I)_I$ is a lacunary family satisfying
$$\|\phi^j_I\|_{\infty} \leq C^{j-1} 2^{-(j-1)N} \quad \text{while} \quad \|\Rest_l\|_{\infty} \leq C^l 2^{- l N} \|f\|_{\infty}.$$
Thus, if $N$ is large enough so that $C\, 2^{-N} < 1$, by {\color{black}letting $l$ go to $\infty$} in (\ref{eq7}), we obtain the desired decomposition (\ref{eq2}) with $\psi_I$ given by
$$\psi_I (x) := \sum_{l=1}^{\infty} \phi^l_I(x).$$
Strictly speaking, the families $(\phi^l_I)_I$ are naturally associated to intervals of length one not $2^{-N}$, but since $N$ is a fixed universal constant, it is not difficult to see that they satisfy 
the estimates (\ref{eq0}) as well, at the expense of loosing a harmless constant of type $2^{1000 N}$. This completes the proof of Proposition \ref{descpr}.

\section{Proof of Theorem \ref{discreteth}}
\label{sec:proof-discretized-main-thm}
Recall that our goal now is to prove that
\begin{equation}\label{discrete-1}
\| \sum_{I\in \I} \langle f, \phi^1_I\rangle \phi^2_I \|_{L^p(L^Q)} \lesssim \Big\| \big( \sum_{I\in\I} \frac{ |\langle f, \phi^1_I\rangle|^2  } {|I|} \one_I \big)^{1/2} \Big\|_{L^p(L^Q)}.
\end{equation}
for every $0<p<\infty$ and every $n$-tuple $Q$ of positive real numbers. Also, $\I$ is a fixed finite collection of dyadic intervals. Of course, the implicit constant in (\ref{discrete-1}) is meant to be independent of the cardinality of $\I$. 
We also denote by $\overline{\I}$ the collection of all dyadic intervals
$J$ having the property that there exists $I\in \I$ so that $I\subseteq J$ and satisfying $|J|\leq 2^M$ for some large fixed positive integer $M$. Sometimes, we refer to the intervals in $\overline{\I}$ as being
the {\it relevant} dyadic intervals.

Let now $E \subseteq \R$ be a measurable subset. To prove (\ref{discrete-1}) it is necessary to prove a more careful version of it given by
\begin{equation}\label{discrete-2}
\| (\sum_{I\in \I} \langle f, \phi^1_I\rangle \phi^2_I ) \one_E \|_{L^p(L^Q)} \lesssim \Big\| \big( \sum_{I\in\I} \frac{ |\langle f, \phi^1_I\rangle|^2  } {|I|} \one_I \big)^{1/2} \Big\|_{L^p(L^Q)} \cdot \big( \size_{\I} \one_E \big)^{1/p-\epsilon}
\end{equation}
where $\epsilon > 0$ is arbitrarily small while
\begin{equation}\label{sss}
\size_{\I} \one_E := \sup_{I\in \overline{\I}} \frac{1}{|I|} \int_{\R} \one_E(x) \big(1+\frac{\dist(x,I) }{|I|}\big)^{-100} d x 
\end{equation}
is essentially the supremum over all $L^1$ averages of $\one_E(x)$ over the intervals of $\overline{\I}$. The reader familiar with our earlier ``helicoidal papers" \cite{BM1}, \cite{BM2} and \cite{BM3} will find our desire to prove (\ref{discrete-2}) natural.

Clearly, (\ref{discrete-2}) implies (\ref{discrete-1}) since one can take $E$ to be the whole real line $\R$.

Using interpolation arguments (see Proposition \ref{prop:interpolation}), it is enough to prove a weaker version of (\ref{discrete-2}), namely
\begin{equation}\label{discrete-3}
\| (\sum_{I\in \I} \langle f, \phi^1_I\rangle \phi^2_I ) \one_E \|_{L^{p,\infty}(L^Q)} \lesssim \big\| \big( \sum_{I\in\I} \frac{ |\langle f, \phi^1_I\rangle|^2  } {|I|} \one_I \big)^{1/2} \big\|_{L^p(L^Q)} \cdot \big( \size_{\I} \one_E \big)^{1/p-\epsilon}.
\end{equation}
Such interpolation arguments will in fact be freely used throughout the section, until the end of it, when they will be proved in detail.

Let us denote by $\P(n)$ the statement which says that (\ref{discrete-3}) holds in full generality, for $0<p<\infty$ and $Q$ $n$-tuple of positive real numbers.
We will prove $\P(n)$ by induction for every $n\geq 0$.

\subsection{Proof of $\P(0)$}

This is the scalar case which now reads as
\begin{equation}\label{discrete-scalar}
\| (\sum_{I\in \I} \langle f, \phi^1_I\rangle \phi^2_I ) \one_E \|_{L^{p,\infty}} \lesssim \Big\| \big( \sum_{I\in\I} \frac{ |\langle f, \phi^1_I\rangle|^2  } {|I|} \one_I \big)^{1/2} \Big\|_{L^p}
\cdot \left( \size_{\I} \one_E \right)^{1/p-\epsilon}.
\end{equation}
Let $s$ be any positive real number with the property $s\leq \min(1,p)$. To estimate the left hand side of (\ref{discrete-scalar}) we dualize the expression through $L^s$, as explained
in \cite{BM2}. Given also the scale invariance of the inequality, this amounts to prove that for every $F\subseteq \R$ measurable set with $|F|=1$, there exists a subset of it $\widetilde{F}\subseteq F$
with $|\widetilde{F}| > 1/2$ so that
\begin{equation}\label{discrete-4}
\| (\sum_{I\in \I} \langle f, \phi^1_I\rangle \phi^2_I ) \one_E  \one_{\widetilde{F}}\|_{L^s} \lesssim \Big\| \big( \sum_{I\in\I} \frac{ |\langle f, \phi^1_I\rangle|^2  } {|I|} \one_I \big)^{1/2} \Big\|_{L^p} \cdot \left( \size_{\I} \one_E \right)^{1/p-\epsilon}.
\end{equation}
To construct the subset $\widetilde{F}$, we start by defining an exceptional set $\Omega$ as follows. 

First, for every integer $k\geq 0$ we define
$$\Omega_k : = \{ x : S f(x) > C 2^{10k/p} \| S f\|_p \}.$$
Here, and from now on, by $S f (x)$ we mean the``discrete" Littlewood-Paley square function given by
\begin{equation}
\label{def:square-func-coll-int}
S f (x) := \big( \sum_{I\in\I} \frac{ |\langle f, \phi^1_I\rangle|^2  } {|I|} \one_I  (x) \big)^{1/2}.
\end{equation}
When we need to emphasize that the square function above depends on the collection $\I$, we write $ S_{\I}$.

It is not difficult to see that
$$|\Omega_k| \leq \frac{1}{2^{10k}} \frac{1}{C^p}.$$
After that we set
$$\widetilde{\Omega_k} : = \{ x : M(\one_{\Omega_k})(x) >  1/ 2^k \}$$
where $M$ is the Hardy-Littlewood maximal operator, and finally
$$\Omega : = \bigcup_{k=0}^{\infty} \widetilde{\Omega_k}.$$
Clearly, 
$$|\widetilde{\Omega_k}| \leq \widetilde{C} \, 2^k |\Omega_k| \leq \widetilde{C}\, 2^k \frac{1}{2^{10k}} \frac{1}{C^p} = \frac{\widetilde{C}}{C^p} \frac{1}{2^{9k}}$$
and in particular this implies that $|\Omega| < 1/10$ if $C$ is a large enough constant\footnote{The constant $\widetilde{C}$  is the boundedness constant of $M : L^1 \rightarrow L^{1,\infty}$.}.

In the end we set $\widetilde{F} : = F\setminus \Omega$ which is a major subset of $F$, in the sense that it satisfies $|\widetilde{F}| \sim 1$. Now, using a result from \cite{MPTT}, we decompose the functions $\phi^2_I$ as
\begin{equation}\label{desco}
\phi^2_I = \sum_{\ell=0}^{\infty} 2^{- M \, \ell} \phi^2_{I, \ell},
\end{equation}
where $M$ is arbitrarily large and for each $\ell \geq 0$, $(\phi^2_{I, \ell})_I$ is still a lacunary family with the additional property that
$$\supp\, \phi^2_{I, \ell} \subseteq 2^\ell I.$$
In particular, one can estimate the left hand side of (\ref{discrete-4}) by
\begin{equation}\label{discrete-5}
\| (\sum_{I\in \I} \langle f, \phi^1_I\rangle \phi^2_I ) \one_E  \one_{\widetilde{F}}\|_s^s \lesssim \sum_{\ell=0}^{\infty} 2^{-Ms \, \ell }
\| (\sum_{I\in \I} \langle f, \phi^1_I\rangle \phi^2_{I,\ell} ) \one_E  \one_{\widetilde{F}}\|_s ^s. 
\end{equation}
The right hand side of (\ref{discrete-5}) can be also rewritten as 
$$\sum_{\ell=0}^{\infty} 2^{-Ms\ell/2} \| (\sum_{I\in \I} \langle f, \phi^1_I\rangle \widetilde{\phi}^2_{I,\ell} ) \one_E  \one_{\widetilde{F}}\|_s ^s$$
where $\widetilde{\phi}^2_{I,\ell} : = 2^{-M\, \ell/2} \phi^2_{I,\ell}$. We will see in what follows that for each $\ell\geq 0$ one has
\begin{equation}\label{discrete-6}
\| (\sum_{I\in \I} \langle f, \phi^1_I\rangle \widetilde{\phi}^2_{I,\ell} ) \one_E  \one_{\widetilde{F}}\|_s ^s \lesssim 2^{L \ell}
\| S f(x) \|_p^s \cdot \left( \size_{\I}\one_E \right)^{(1/p-\epsilon) s}
\end{equation}
where $L$ is some constant depending on $s$ and $p$. However, because of the large constant $M$ in \eqref{discrete-5}, this will be enough to complete our proof.
We will prove (\ref{discrete-6}) in detail in the main case when $\ell=0$ and then we will explain how to modify the argument to obtain \eqref{discrete-6} in general.

In other words, the goal for us now is to prove that
\begin{equation}\label{discrete-7}
\| (\sum_{I\in \I} \langle f, \phi^1_I\rangle \widetilde{\phi}^2_{I,0} ) \one_E  \one_{\widetilde{F}}\|_s ^s \lesssim 
\| S f(x) \|_p^s \cdot \left( \size_{\I} \one_E \right)^{(s/p-\epsilon)}.
\end{equation}
Recall that now, since 
$$\supp\, \widetilde{\phi}^2_{I, 0} \subseteq I$$
one must have $I \cap \Omega^c \neq \emptyset$ which in particular implies that $I\cap \Omega_0^c \neq \emptyset$. From the definition of $\Omega_0$, one can see that this set
admits a natural decomposition as a disjoint union of maximal dyadic intervals denoted by $I_{max}$. In particular, our dyadic intervals $I$ have the property that they are either disjoint from all
these $I_{max}$, or they contain strictly at least one of them. In either case, it is not difficult to see that one has the pointwise estimate
\begin{equation}\label{point}
\big( \sum_{I\in\I :  I\cap \Omega_0^c \neq \emptyset   } \frac{ |\langle f, \phi^1_I\rangle|^2  } {|I|} \one_I (x)  \big)^{1/2} \leq \widetilde{C} \| S  f \|_p
\end{equation}
where $\widetilde{C}$ is a universal constant.
To prove (\ref{discrete-7}) we will combine two stopping time arguments, one performed with the help of averages of the type
\begin{equation}\label{ave1}
\frac{1}{|I_0|^{1/p}}
\Big\| \big( \sum_{I\subseteq I_0} \frac{ |\langle f, \phi^1_I\rangle|^2  } {|I|} \one_I \big)^{1/2} \Big\|_{L^p}
\end{equation}
and the other with the help of averages of type
\begin{equation}\label{ave2}
\frac{1}{|I_0|} \int_{\R} \one_{E\cap\widetilde{F}} (x) \big(1+ \frac{ \dist(x, I_0) } {|I_0|}\big)^{-100} d x. 
\end{equation}
The latter will be denoted from now on $\ave^1_{I_0} (\one_{E\cap\widetilde{F}})$. Clearly, because of the pointwise bound (\ref{point}), averages such as the ones in (\ref{ave1})
cannot be larger than $\widetilde{C}\, \| S f \|_p$, while averages of type (\ref{ave2}) cannot be larger than $\size_{\I}(\one_{E\cap \widetilde{F}})$.

We describe now in detail the first stopping time.

We start by selecting maximal dyadic intervals $I_0\in \overline{\I}$ with the property that $I_0\cap \Omega^c_0\neq \emptyset$ and so that
\begin{equation}\label{eq8}
\frac{1}{|I_0|^{1/p}}
\Big\| \big( \sum_{I\subseteq I_0} \frac{ |\langle f, \phi^1_I\rangle|^2  } {|I|} \one_I \big)^{1/2} \Big\|_{L^p} \geq \frac{ \widetilde{C}}{2} \| S f\|_p.
\end{equation}
Of course, as pointed out before, we implicitly assume that all the intervals $I$ that participate in the summation above have the property $ I \cap \Omega_0^c \neq \emptyset$. It is also important to observe that these selected intervals $I_0$ are all disjoint, as a consequence of their maximality. Then, we ignore all the relevant dyadic intervals that lie inside one of these selected intervals and consider only those that are left. They are either disjoint
from the selected ones, or they contain at least one of the selected ones.

After this, among those that are left, we pick those maximal ones, still denoted by $I_0$ for which
\begin{equation}\label{eq9}
\frac{1}{|I_0|^{1/p}}
\Big\| \big( \sum_{I\subseteq I_0} \frac{ |\langle f, \phi^1_I\rangle|^2  } {|I|} \one_I \big)^{1/2} \Big\|_{L^p} \geq \frac{ \widetilde{C}}{2^2} \| S f\|_p
\end{equation}
and so forth. The maximal intervals selected at the first step are collected in $\I^{(1)}_1$, those selected at the second step are collected in $\I^{(1)}_2$ and so on, obtaining the collections
$(\I^{(1)}_{n_1})_{n_1}$. Clearly, there are only finitely many such steps, since our initial collection of intervals was finite.

After that, independently, we perform a similar stopping time, but one that involves the averages $\ave_{I_0}^1(\one_{E\cap\widetilde{F}})$ instead. We start by selecting those maximal intervals $I_0$ for which
\[
\ave_{I_0}^1(\one_{E\cap\widetilde{F}}) > \frac{1}{2} \size_{\I}(\one_{E\cap\widetilde{F}})
\]
then, among those that are left (more specifically, those that are not inside any of the previously selected $I_0$) we pick again those maximal $I_0$ for which
$$\ave_{I_0}^1(\one_{E\cap\widetilde{F}}) > \frac{1}{2^2} \size_{\I}(\one_{E\cap\widetilde{F}})$$
and so on, exactly as before. In this way, one obtains a sequence of collections of maximal dyadic intervals $I_0$ denoted by $(\I^{(2)}_{n_2})_{n_2}$.

In the end, we combine them to be able to estimate (\ref{discrete-7}). One can write
\begin{align}
\label{eq10}
&\| (\sum_{I\in\I, I\cap\Omega^c_0\neq\emptyset} \langle f, \phi^1_I \rangle \widetilde{\phi}^2_{I,0}) \one_E \one_{\widetilde{F}} \|_s^s  \leq \sum_{n_1, n_2} \sum_{I_1\in \I^{(1)}_{n_1}, I_2\in\I^{(2)}_{n_2}}
\| (\sum_{I \in  \I^{(1)}_{n_1}(I_1) \cap \I^{(2)}_{n_2}(I_2) } \langle f, \phi^1_I \rangle \widetilde{\phi}^2_{I,0}) \one_E \one_{\widetilde{F}} \|_s^s
\end{align}
where $\I^{(1)}_{n_1}(I_1)$ contains all the relevant dyadic intervals $I$ with the property that $I\subseteq I_1$ but such that $I$ is not contained in any of the previously selected intervals in $\I^{(1)}_l$
for $0\leq l\leq n_1-1$, and similarly for $\I^{(2)}_{n_2}(I_2)$. Clearly, any interval $I$ participating in the summation (\ref{eq10}) must satisfy $I\subseteq I_1\cap I_2$. Now, for every $I_1, I_2$ as before,
the corresponding $L^s$ quasi-norm in (\ref{eq10}) can be estimated by
\begin{equation}\label{eq11}
\| (\sum_{I\subseteq I_1\cap I_2   } \langle f, \phi^1_I \rangle \widetilde{\phi}^2_{I,0}) \one_E \one_{\widetilde{F}} \|_1 \cdot |E\cap \widetilde{F} \cap I_1 \cap I_2 |^{\frac{1-s}{s}}
\end{equation}
by using H\"{o}lder, since $s\leq 1$. The $L^1$ norm in (\ref{eq11}) can be dualized and estimated by
$$\sum_{I\subseteq I_1\cap I_2}
\langle f, \phi^1_I \rangle \langle \one_{E\cap\widetilde{F}} \, g, \widetilde{\phi}^2_{I,0} \rangle
$$
for some function $g$ with the property $\|g\|_{\infty}=1$. Using Cauchy-Schwartz this can be further estimated by
{\fontsize{10}{10}$$
\frac{1}{ |I_1\cap I_2 |^{1/2}}
\big\| \big( \sum_{I\subseteq I_1\cap I_2} \frac{ |\langle f, \phi^1_I\rangle|^2  } {|I|} \one_I \big)^{1/2} \big\|_2 \cdot  \frac{1}{ |I_1\cap I_2 |^{1/2}}
\big\| \big( \sum_{I\subseteq I_1\cap I_2} \frac{ |\langle     \one_{E\cap\widetilde{F}} \, g, \widetilde{\phi}^2_{I,0}      \rangle|^2  } {|I|} \one_I \big)^{1/2} \big\|_2 \cdot |I_1\cap I_2|.
$$}
Using John-Nirenberg now twice (see Theorem $2.10$ in \cite{cw} for this robust, discrete, variant of it) together with the standard local estimate of weak-$L^1$ averages (which can be found in Lemma 2.16 of \cite{cw} for instance), this can be further majorized by 
\begin{equation}
\Big(   \sup_{J_1 \subseteq I_1\cap I_2}   \frac{1}{|J_1|^{1/p}}
\big\| \big( \sum_{I\subseteq J_1} \frac{ |\langle f, \phi^1_I\rangle|^2  } {|I|} \one_I \big)^{1/2} \big\|_p \Big) \cdot \big(  \sup_{J_2 \subseteq I_1\cap I_2} \ave^1_{J_2} (\one_{E\cap\widetilde{F}}) \big) \cdot |I_1\cap I_2|.
\end{equation}

If one raises these estimates to the power $s$, as required by (\ref{eq10}), one can see that the corresponding expression there is smaller than
\begin{align}
&\Big(   \sup_{J_1 \subseteq I_1\cap I_2}   \frac{1}{|J_1|^{1/p}}
\big\| \big( \sum_{I\subseteq J_1} \frac{ |\langle f, \phi^1_I\rangle|^2  } {|I|} \one_I \big)^{1/2} \big\|_p   \Big)^s \cdot \big(  \sup_{J_2\subseteq I_1\cap I_2} \ave^1_{J_2} (\one_{E\cap\widetilde{F}})   \big)^s \cdot |I_1\cap I_2|^s \\
&\quad \cdot \big(    \ave^1_{I_1\cap I_2} (\one_{E\cap \widetilde{F}})  \big)^{1-s} \cdot |I_1\cap I_2|^{1-s}, \nonumber
\end{align}
which is smaller still than
\begin{equation}
\Big(   \sup_{J_1 \subseteq I_1\cap I_2}   \frac{1}{|J_1|^{1/p}}
\big\| \big( \sum_{I\subseteq J_1} \frac{ |\langle f, \phi^1_I\rangle|^2  } {|I|} \one_I \big)^{1/2} \big\|_p  \Big)^s \cdot \big(  \sup_{J_2\subseteq I_1\cap I_2} \ave^1_{J_2} (\one_{E\cap\widetilde{F}})   \big) \cdot |I_1\cap I_2|.
\end{equation}

Using these estimates in (\ref{eq10}) the expression there can be estimated further by
\begin{equation}\label{eq12}
\sum_{n_1, n_2} \sum_{I_1\in\I^{(1)}_{n_1}, \, I_2\in\I^{(2)}_{n_2}} 2^{-n_1 s} 2^{-n_2} |I_1\cap I_2|.
\end{equation}
On the other hand the expression 
$$\sum_{I_1\in\I^{(1)}_{n_1}, \, I_2\in\I^{(2)}_{n_2}} |I_1\cap I_2|$$
is smaller than
$$\sum_{I_1\in\I^{(1)}_{n_1}} |I_1| \lesssim 2^{n_1 p} \|S f \|_p^p$$
and also smaller than
$$\sum_{I_2\in\I^{(2)}_{n_2}} |I_2| \lesssim 2^{n_2}$$
given that $|\widetilde{F}| \sim 1$. This implies that
$$\sum_{I_1\in\I^{(1)}_{n_1}, \, I_2\in\I^{(2)}_{n_2}} |I_1\cap I_2| \lesssim 2^{n_1 p \theta_1}  \|S f\|_p^{p\theta_1}  2^{n_2 \theta_2}$$
for every $0\leq \theta_1, \theta_2 \leq 1$ so that $\theta_1+\theta_2 = 1$. Using this in (\ref{eq12}) one can majorize that expression by
\begin{equation}\label{eq13}
\sum_{n_1, n_2} 2^{-n_1 (s - p \theta_1)} 2^{ -n_2 (1 - \theta_2)} \| S f \|_p^{p \theta_1}.
\end{equation}
But now, we recall that $2^{- n_1} \lesssim \| S f\|_p$ while $2^{- n_2} \lesssim \size_{\I} \one_E$ and in particular this means that (\ref{eq13}) is smaller than
$$\| S f \|_p^{s-p \theta_1} \cdot \|S f \|_p^{p \theta_1} \cdot \big(\size_{\I} \one_E \big)^{1 - \theta_2}$$
provided that $\theta_1 < s/p$ which is the condition that guarantees the convergence of the geometric series over $n_1$. If $\theta_1$ is taken very close to $s/p$, this gives an upper bound
of type
$$\| S f \|_p^s \cdot \big( \size_{\I}\one_E \big)^{s/p - \epsilon}$$
as desired in (\ref{discrete-7}). 

To prove (\ref{discrete-6}) for arbitrary $\ell > 0$ one proceeds similarly. The observation now is that since $\supp \, \widetilde{\phi}^2_{I, \ell} \subseteq 2^\ell I$ one must have
$$2^\ell I \cap \Omega^c \neq \emptyset$$
and it is not difficult to see that this implies that
$$I \cap \Omega^c_\ell \neq \emptyset.$$
Indeed, if this was not true, then $I\subseteq \Omega_\ell$, which means that $2^\ell I \subseteq \widetilde{\Omega}_\ell \subseteq \Omega$, a contradiction.

Now one simply repeats the argument before. One difference is that the first $L^p$ averages of the square function can be as large as $C 2^{10 \ell /p} \| S f\|_p$, a bound which is responsible for the positive constant $L$ in (\ref{discrete-6}). Another difference is in the estimate (\ref{eq11}), whose analogue now  contains a factor of type
$$|E \cap \widetilde{F} \cap 2^\ell (I_1 \cap I_2)|^{\frac{1-s}{s}}.$$
However, the small constant $2^{-M \ell /2}$ in the definition of $\widetilde{\phi}^2_{I, \ell}$ gets multiplied by it, and this allows one to write
$$2^{-M \ell /2} |E \cap \widetilde{F} \cap 2^\ell (I_1 \cap I_2)|^{\frac{1-s}{s}}\lesssim \big( \int_{\R} \one_{E\cap \widetilde{F}} (x) \big(  1+ \frac{\dist (x, I_1\cap I_2)  }{|I_1 \cap I_2| }   \big)^{-100}  d x \big)^{\frac{1-s}{s}}$$
and everything continues as before, if $M$ is large enough. This completes the proof of $\P(0)$.

\subsection{Proof of $\P(n-1)$ implies $\P(n)$.}
\label{sec:inductive-section}
Recall that what we need to prove now is the estimate
\begin{equation}\label{eq14}
\| (\sum_{I\in \I} \langle f, \phi^1_I\rangle \phi^2_I ) \one_E \|_{L^{p,\infty}(L^Q)} \lesssim \Big\| \big( \sum_{I\in\I} \frac{ |\langle f, \phi^1_I\rangle|^2  } {|I|} \one_I \big)^{1/2} \Big\|_{L^p(L^Q)} \cdot \big( \size_{\I} \one_E \big)^{1/p-\epsilon} 
\end{equation}
for every $0 < p < \infty$ and $Q$ an $n$-tuple of positive real numbers, assuming that even the stronger version of it, namely (\ref{discrete-2}), holds true for $(n-1)$-tuples $Q$. Again, here we are implicitly assuming that
the proof of the strong $L^p(L^Q)$ estimate in (\ref{eq14}) will follow by standard interpolation arguments, which we will describe later on, as promised.

Define $q_{j_0}:= \min\limits_{1\leq j \leq n} q_j$ and let $s$ be any positive real number so that $s\leq \min (1, p, q_{j_0})$. Then, one can dualize the  weak-$L^p$ quasi-norm on the left hand side of
(\ref{eq14}) through $L^s$, as explained in \cite{BM2}. As before, this amounts to prove that for every $F\subseteq\R$ measurable set with $|F|=1$, there exists a subset $\widetilde{F} \subseteq F$ with
$|\widetilde{F}| > 1/2$ so that
\begin{equation}\label{eq15}
\| \|(\sum_{I\in \I} \langle f, \phi^1_I\rangle \phi^2_I ) \one_E \one_{\widetilde{F}}  \|_Q  \|_s \lesssim  \|S f \|_{L^p(L^Q)}
\cdot \left( \size_{\I} \one_E  \right)^{1/p-\epsilon}.
\end{equation}
To construct $\widetilde{F}$, one first constructs an exceptional set $\Omega$, as in the scalar case, with the only difference that the corresponding $\Omega_k$ is given now by
$$\Omega_k : = \{ x : \|S f(x) \|_Q > C \, 2^{10 k / p} \| \| S f \|_Q \|_p \}.$$

After that, one defines $\widetilde{F} := F \setminus \Omega$ exactly as before, which is clearly  a major subset of $F$, in the sense that it has a comparable measure. Then, one uses again the decomposition (\ref{desco})
to reduce matters to proving the analogue of (\ref{discrete-6}) which is now given by
\begin{equation}\label{eq16}
\|  \| (\sum_{I\in \I} \langle f, \phi^1_I\rangle \widetilde{\phi}^2_{I,k} ) \one_E  \one_{\widetilde{F}} \|_Q \|_s ^s \lesssim 2^{Lk}
\| \| S f \|_Q \|_p^s \cdot \left( \size_{\I} \one_E \right)^{(1/p-\epsilon) s}.
\end{equation}
Recall from \cite{BM2} that $s\leq \min (1, p, q_{j_0})$ implies that the expression on the left hand side of (\ref{eq16}) is sub-additive now. As before, we will describe the proof of (\ref{eq16}) in the main case $k=0$,
the changes in the general case being similar to the ones in the scalar case. We therefore want to show that
\begin{equation}\label{eq17}
\|  \| (\sum_{I\in \I} \langle f, \phi^1_I\rangle \widetilde{\phi}^2_{I,0} ) \one_E  \one_{\widetilde{F}} \|_Q \|_s ^s \lesssim 
\| \| S f \|_Q \|_p^s \cdot \left( \size_{\I} \one_E \right)^{(1/p-\epsilon) s}.
\end{equation}
To estimate the left hand side of (\ref{eq17}) we combine as before, two stopping times. The first one, selects iteratively maximal dyadic intervals $I_0$ for which one has
\begin{equation}\label{aveQ}
\frac{1}{|I_0|^{1/p}} \,
\Big\| \big\| \big( \sum_{I\subseteq I_0} \frac{ |\langle f, \phi^1_I\rangle|^2  } {|I|} \one_I \big)^{1/2} \big\|_Q \Big\|_{L^p} >  \frac{\widetilde{C}}{2^l} \, \| S f \|_{L^p(L^Q)}
\end{equation}
for various $l\geq 0$, while the second is identical to the one used in the scalar case. This allows us to estimate the left hand side of (\ref{eq17}) by

\begin{equation}\label{eq18}
\sum_{n_1, n_2} \sum_{I_1\in \I^{(1)}_{n_1}, I_2\in\I^{(2)}_{n_2}}
\|  \|  (\sum_{I\in  \I^{(1)}_{n_1}(I_1) \cap \I^{(2)}_{n_2}(I_2)     } \langle f, \phi^1_I \rangle \widetilde{\phi}^2_{I,0}) \one_E \one_{\widetilde{F}}  \|_Q \|_s^s.
\end{equation}
Fix now $I_1$ and $I_2$ and consider the corresponding term on the right hand side of (\ref{eq18}). Given variables $(s_1, ..., s_n)\in S$ denote by
$\widetilde{s}:= (s_2, ..., s_{n})$ and given $Q= (q_1, ..., q_n)$ denote by $\widetilde{Q} := (q_2, ..., q_{n})$. Using these notations, the expression becomes
{\fontsize{10}{10}\begin{align*}
&\int_{\R} \big\|    (\sum_I  \langle f, \phi^1_I \rangle \widetilde{\phi}^2_{I,0}) \one_E \one_{\widetilde{F}}   \big\|^s_Q (x) d x 
=\int_{\R} \big( \int_{S_1}   \|    (\sum_I  \langle f_{(s_1, \widetilde{s})}   , \phi^1_I \rangle \widetilde{\phi}^2_{I,0}) \one_E \one_{\widetilde{F}}  \|_{L^{\widetilde{Q}}_{\widetilde{s}}}^{q_1} (x) d s_1       \big)^{s/q_1} d x  \\
=&\int_{\R} \big( \int_{S_1}  \|    (\sum_I  \langle f_{(s_1, \widetilde{s})}   , \phi^1_I \rangle \widetilde{\phi}^2_{I,0}) \one_E \one_{\widetilde{F}} \|_{L^{\widetilde{Q}}_{\widetilde{s}}}^{q_1} (x) d s_1    \big)^{s/q_1}   \one_E (x) \one_{\widetilde{F}}  (x)  d x .
\end{align*}}

Since $s/q_1\leq 1$ one can apply H\"{o}lder and estimate the above expression by
\begin{equation}\label{eq19}
\big[   \int_{\R} \int_{S_1}  \|    (\sum_I  \langle f_{(s_1, \widetilde{s})}   , \phi^1_I \rangle \widetilde{\phi}^2_{I,0})  \one_E \one_{\widetilde{F}}    \|_{L^{\widetilde{Q}}_{\widetilde{s}}}^{q_1} (x) d s_1     d x    \big]^{s/q_1} \cdot
|E \cap \widetilde{F} \cap I_1 \cap I_2 |^{(1- s/q_1)} 
\end{equation}
using also the fact that all the intervals $I$ are now inside $I_1\cap I_2$.
Then, one can use Fubini and integrate first with respect to the $x$ variable in (\ref{eq19}). This allows one to use the induction hypothesis locally (i.e. with respect to the collection  $\I^{(1)}_{n_1}(I_1) \cap \I^{(2)}_{n_2}(I_2) $ )
in the case $p=q_1$, and estimate (\ref{eq19}) by
\begin{equation}\label{eq20}
\big(  \int_{\R} \int_{S_1} \| S f (x) \|_{\widetilde{Q} }^{q_1}d s_1  d x    \big)^{\frac{s}{q_1}} 
\cdot
\big(   \size_{ \I^{(1)}_{n_1}(I_1) \cap \I^{(2)}_{n_2}(I_2)     }   (\one_{E \cap \widetilde{F}})       \big)^{(\frac{s}{q_1} - \epsilon)}
\cdot
|E \cap \widetilde{F} \cap I_1 \cap I_2 |^{(1- \frac{s}{q_1})} .
\end{equation} 
We emphasize that in (\ref{eq20}) the implicit sum in the definition of the square function $S f(x)$ runs over the intervals $I$ inside the local collection  $\I^{(1)}_{n_1} (I_1) \cap \I^{(2)}_{n_2} (I_2)$. 

It is then not difficult to see that the last expression in (\ref{eq20}) can be rewritten and majorized by 
{\fontsize{11}{10}\begin{align}\label{eq21}
&\Big(  \frac{1}{| I_1 \cap I_2 |^{1/q_1}} \,\big\| \big\| \big( \sum_{I\subseteq I_1\cap I_2} \frac{ |\langle f, \phi^1_I\rangle|^2  } {|I|} \one_I \big)^{1/2} \big\|_Q \big\|_{L^{q_1}}   \Big)^s \cdot \big(     \size_{   \I^{(1)}_{n_1}(I_1) \cap \I^{(2)}_{n_2}(I_2)     }   (\one_{E \cap \widetilde{F}})   \big)^{1-\epsilon} \cdot |I_1\cap I_2| . 
\end{align}}

Using once again the John-Nirenberg inequality from \cite{cw} (which works equally well in our multiple vector-valued setting), we find that (\ref{eq21}) is smaller than
\begin{equation}\label{eq22}
\sup_{J\subseteq I_1\cap I_2} \Big(     \frac{1}{| J |^{1/p}}
\Big\| \big\| \big( \sum_{I\subseteq J } \frac{ |\langle f, \phi^1_I\rangle|^2  } {|I|} \one_I \big)^{1/2} \big\|_Q \Big\|_{L^{p}}          \Big)^s \cdot 
\big(  \size_{   \I^{(1)}_{n_1}(I_1) \cap \I^{(2)}_{n_2}(I_2)     }   (\one_{E \cap \widetilde{F}})          \big)^{1-\epsilon} \cdot |I_1\cap I_2| .
\end{equation}
Using these, we can come back to (\ref{eq18}) and majorize that expression by
$$\sum_{n_1, n_2} 2^{-n_1 s} 2^{- n_2 (1-\epsilon)} \sum_{I_1\in \I^{(1)}_{n_1},  I_2\in \I^{(2)}_{n_2}} |I_1 \cap I_2|.$$
As before, one can estimate 
$$    \sum_{I_1\in \I^{(1)}_{n_1},  I_2\in \I^{(2)}_{n_2}} |I_1 \cap I_2|       $$ 
in two distinct ways, by taking advantage of the stopping time decompositions performed earlier.

First, we can estimate it by $2^{n_1 p} \| \|S f \|_Q \|_p^p$ and second, by $2^{n_2}$ given that $|\widetilde{F}| \sim 1$. In particular, this allows one to estimate the whole expression by
$$\sum_{n_1, n_2} 2^{- n_1 (s - p \theta_1)} 2^{- n_2 (1- \epsilon - \theta_2)} \| S f \|_{L^p(L^Q)}^{p\theta_1}$$
as in the scalar case, for every $ 0\leq \theta_1, \theta_2 \leq 1$ with   $\theta_1 + \theta_2 = 1$. Then, if one chooses $\theta_1 < s/p$ but very close to it, this double sum becomes smaller than
$$\| S f \|^{s - p \theta_1}_{L^p(L^Q)} \cdot \left(  \size_{\I} \one_E  \right)^{s/p - \epsilon} \cdot \| S f \|_{L^p(L^Q)}^{p\theta_1}$$
as desired. And this completes our proof.

The only thing left is the interpolation argument that we used implicitly several times.

\subsection{Interpolation}
Our interpolation result is somewhat unusual, in the sense that the collection $\ii I$ of dyadic intervals is as important as the operator it defines, the square function associated to it from \eqref{def:square-func-coll-int}. The result and its proof generalize straight away to collections of cubes in $\R^d$, and to arbitrary measures.

\begin{proposition}
\label{prop:interpolation}
Consider $0 < p_1 < p < p_2 < \infty$ and let $\tilde {\mathcal{I}}$ be a collection of dyadic intervals. Assume that, for any subcollection $\ii I \subseteq \tilde{\mathcal{I}}$ of dyadic intervals and any $L^Q$-valued Schwartz function $f$ on $\rr R$, we have for $j= 1, 2$,
\begin{equation}\label{inter1}
\| (\sum_{I \in \ii I} \langle f, \phi^1_I\rangle \phi^2_I ) \one_E \|_{L^{p_j,\infty}(L^Q)} \leq K_j \,  \big\| \big( \sum_{I\in \ii I} \frac{ |\langle f, \phi^1_I\rangle|^2  } {|I|} \one_I \big)^{1/2} \big\|_{L^{p_j}(L^Q)}
\end{equation}
with the constants $K_j$ independent on $\ii I$. Then for any $\I \subseteq \tilde{\I}$ we have the strong bound
\begin{equation}\label{inter2}
\| (\sum_{I \in \ii I} \langle f, \phi^1_I\rangle \phi^2_I ) \one_E\|_{L^p(L^Q)} \leq K \, \big\| \big( \sum_{I\in \ii I} \frac{ |\langle f, \phi^1_I\rangle|^2  } {|I|} \one_I \big)^{1/2} \big\|_{L^p(L^Q)},
\end{equation}
where $K \lesssim \big( K_1^{p_1} +K_2^{p_2} \big)^\frac{1}{p} $.

\begin{observation}
As mentioned before, the interpolation result in Proposition \ref{prop:interpolation} can be stated in a more general setting, as the interested reader can verify. Our choice of presentation is motivated by the fact that in the present paper we need precisely the form presented above. The constants $K_j$ in \eqref{inter1} do not depend on any of the subcollections $\ii I$ of intervals, but they could (and in most applications they do) depend on the general collection $\tilde{\I}$ and on the set $E$ appearing on the left hand side of \eqref{inter1}; as a consequence, $K$ is not dependent upon any of the subcollections $\ii I$, but could depend on $\tilde{\I}$ and on the set $E$.

We use the interpolation result above in order to deduce \eqref{discrete-2}  from \eqref{discrete-3}; notice that in that case 
\[
K_j= \big( \size_{\tilde{\I}} \one_E   \big)^{\frac{1}{p_j}-\epsilon} \quad \text{for    } j=1, 2.
\]
Hence \eqref{discrete-2} follows immediately from \eqref{discrete-3} after interpolating carefully in a small neighborhood of the desired index $0 < p < \infty$.

On the other hand, in the proof of the interpolation result we will assume that $E$ is the entire real line since it doesn't play a role in the interpolation argument.
\end{observation}

\begin{proof}[Proof of Proposition \ref{prop:interpolation}]

Let $\ii  I \subseteq \tilde{\I}$ be a subcollection of dyadic intervals, and denote by $F(x)$ the $L^Q$-valued function 
\[
F(x):=\sum_{I \in \ii I} \langle f, \phi_I^1  \rangle \phi_I^2(x).
\]

Our goal is to control $\ds \|F \|_{L^p(L^Q)}$ by $\ds \|S_{\ii I} f \|_{L^p(L^Q)}$, where $S_{\ii I} f$ is the associated square function:
\[
S_{\ii I} f(x):= \big( \sum_{I\in \ii I} \frac{ |\langle f, \phi^1_I\rangle|^2  } {|I|} \one_I \big)^{1/2} .
\]

The proof that we are about to provide will involve a partitioning of the collection $\ii I$ according to level sets of the ``global" square function $S_{\ii I}$. First, for any $\alpha>0$ and any $k \geq 0$, we define
\begin{equation}
\label{eq:def:S:alpha:k)}
S(k, \alpha):= \Big \lbrace x: \|  S_{\ii I}f(x)  \|_{L^Q}>\frac{\alpha}{C^k}   \Big \rbrace,
\end{equation}
where $C$ is a constant that will be determined later. Notice that the sets $S(k, \alpha)$ are nested:
\[
S(0, \alpha) \subseteq S(1, \alpha) \subseteq \ldots S(k, \alpha) \subseteq \ldots.
\]

Each of the sets $S(k, \alpha)$ can be written as a disjoint union of maximal dyadic intervals: 
\[
S(k, \alpha):= \bigcup_{I_{max}^k \in \ii M_k} I_{max}^k, \qquad \forall \, k \geq 0. 
\]
These will be used for the formerly mentioned partition:
\begin{enumerate}[leftmargin=*]
\item[-] the collection $\ii I_0$ will consist of all intervals $I \in \ii I$ that are contained inside some maximal interval $I_{max}^0$:
\[
\ii I_0:= \{ I \in \ii I:  \text{    there exists some    } I_{max}^0 \in \ii M_0 \text{    with    } I \subseteq I_{max}^0 \subset S(0, \alpha)  \}.
\]
\item[-] for any $k \geq 1$, $\ii I_k$ is defined as
{\fontsize{11}{11}\[
\ii I_k:= \{ I \in \ii I:  \exists I_{max}^k \in \ii M_k \text{    with    } I \subseteq I_{max}^k \subset S(k, \alpha) \text{  and  } I \nsubseteq  I_{max}^\ell \text{   for all  } 0 \leq \ell < k  \}.
\]} 
That is, $\I_k$ consists of all the intervals in $\I$ contained in some $I_{max}^k \in \ii M_k$, which were not previously selected in any other $\I_\ell$ with $0 \leq \ell \leq k-1$.
\end{enumerate}

Then we have $\ds \ii I= \bigcup_{k \geq 0} \ii I_k$ and if $F_k$ denotes the $L^Q$-valued function 
\[
F_k(x):=\sum_{I \in \ii I_k} \langle f, \phi_I^1  \rangle \phi_I^2(x),
\]
we have the decomposition $\ds F(x)=\sum_{k \geq 0} F_k(x)$.

Notice that for all $k \geq 0$
\[
\supp ( S_{\ii I_k} f ) \subseteq \bigcup_{I \in \ii I_k} I \subseteq \bigcup_{I_{max}^k \in \ii M_k} I_{max}^k = S(k, \alpha).
\]

For any $k \geq 1$, if $I \in \ii I_k$ and $I_{max}^{k-1} \in \ii M_{k-1}$ are so that $\ds I \cap I_{max}^{k-1} \neq \emptyset$, then necessarily $\ds  I \supsetneq I_{max}^{k-1}$. Given the maximality condition based on which $I_{max}^{k-1}$ was selected in $\ii M_{k-1}$, all intervals $I \in \ii I_k$ intersect $\ds S(k-1, \alpha)^c$ and
\begin{equation}
\label{eq:SF-est-S:k:aplh}
\|S_{\ii I_k} f(x)\|_{L^Q} \leq \|S_{\ii I}f (x)\|_{L^Q} \cdot \one_{S(k-1, \alpha)^c}(x) \leq \frac{\alpha}{C^{k-1}},
\end{equation}
a feature that will be exploited later on.

Since $Q$ is an arbitrary $n$-tuple of positive real numbers, there is no certainty that $\ds \| \cdot \|_{L^Q}$ satisfies the triangle inequality; however, for $s$ small enough (the condition that $\ds s \leq \min_{1 \leq j \leq n} q_j$ suffices), $\ds  \| \cdot \|_{L^Q}^s$ becomes subadditive. As a result, 
\[
\big  \lbrace x : \|\sum_{k \geq 0} F_k(x)     \|_{L^Q}^s > \alpha^s   \big \rbrace  \subseteq \bigcup_{k \geq 0} \big  \lbrace x : \| F_k(x)    \|_{L^Q}^s > \frac{\alpha^s}{2^{k+1}}  \big \rbrace.
\]
Moreover,
\begin{align*}
& \big| \lbrace x: \| \sum_{k \geq 0} F_k(x)  \|_{L^Q}>\alpha  \rbrace  \big|=\big| \lbrace x: \| \sum_{k \geq 0} F_k(x)  \|_{L^Q}^s>\alpha^s  \rbrace  \big| \\
& \leq \sum_{k \geq 0} \big| \lbrace x: \|  F_k(x)  \|_{L^Q}^s>\frac{\alpha^s}{2^{k+1}}  \rbrace  \big|= \sum_{k \geq 0} \big| \lbrace x : \|  F_k(x)  \|_{L^Q}> \frac{\alpha}{2^{(k+1)/s}}  \rbrace  \big|.
\end{align*}

Such an inequality is important because it allows us to estimate $\|F \|_{L^p(L^Q)}$:
\begin{align}\label{eq:interp:subadd}
\|F \|_{L^p(L^Q)}^p & = p \int_{0}^\infty \alpha^{p-1}  \big| \lbrace x: \| \sum_{k \geq 0} F_k(x)  \|_{L^Q}>\alpha  \rbrace  \big| d \alpha \\
& \leq  \sum_{k \geq 0} p \int_{0}^\infty \alpha^{p-1}  \big| \lbrace x : \|  F_k(x)  \|_{L^Q}> \frac{\alpha}{2^{(k+1)/s}}  \rbrace  \big| d \alpha. \nonumber
\end{align}

We note that the functions $F_k$ above depend in fact on the variable $\alpha$ (the collections of intervals $\ii I_k$ are determined by the level sets $S(k, \alpha)$); this is the main difficulty in proving the interpolation result, and which differentiates Proposition \ref{prop:interpolation} from standard interpolation results.

First, we deal with the case corresponding to $k=0$, by invoking the weak-type hypothesis \eqref{inter1} for the collection $\ii I_0$:
\begin{align*}
&  \int_{0}^\infty \alpha^{p-1}  \big| \lbrace x : \|  F_0(x)  \|_{L^Q}> \frac{\alpha}{2^{1/s}}  \rbrace  \big| d \alpha \\
& \leq \int_0^\infty \alpha^{p-1} \Big(  \frac{\alpha}{2^{1/s}} \Big)^{-p_1} \, K_1^{p_1} \big\| \big( \sum_{I \in \ii I_0} \frac{|\langle f, \phi_I^1 \rangle|^2}{|I|} \cdot \one_I  \big)^\frac{1}{2} \big\|_{L^{p_1}(L^Q)}^{p_1} \, d \alpha \\
& \leq 2^{p_1/s} K_1^{p_1} \int_0^\infty \alpha^{p-p_1-1} \int_{S(0, \alpha)} \| S_{\ii I_0} f(x)   \|_{L^Q}^{p_1} \, d x d \alpha. 
\end{align*}

The term above can be bounded by an expression involving only the ``global" square function $S_{\I}$, which doesn't depend on the subcollection $\I_0$ nor on $\alpha$, given by
\[
2^{p_1/s} K_1^{p_1} \int_0^\infty \alpha^{p-p_1-1} \int_{\lbrace \|S_{\ii I} f_{L^Q}\|_{L^Q} > \alpha   \rbrace}  \|S_{\ii I} f(x) \|_{L^Q}^{p_1} \, d x d \alpha.
\]
Now we apply the usual trick which consists in changing the order of integration, obtaining in this way
\begin{align*}
&  \int_{0}^\infty \alpha^{p-1}  \big| \lbrace x : \|  F_0(x)  \|_{L^Q}> \frac{\alpha}{2^{1/s}}  \rbrace  \big| d \alpha \\ & \leq 2^\frac{p_1}{s} K_1^{p_1} \int_{\rr R}  \|S_{\ii I} f(x) \|_{L^Q}^{p_1} \int_0^{\|S_{\ii I} f(x) \|_{L^Q}} \alpha^{p-p_1-1} d \alpha = \frac{ 2^\frac{p_1}{s} K_1^{p_1}}{p-p_1} \|S_{\ii I} f \|_{L^p(L^Q)}^p.
\end{align*}

Next we deal with a generic term involving $F_k$ for some $k  \geq 1$; we use the assumption \eqref{inter1} applied to the collection $\ii I_k$:
{\fontsize{11}{11}\begin{align*}
&  \int_{0}^\infty \alpha^{p-1}  \big| \lbrace x : \|  F_k(x)  \|_{L^Q}> \frac{\alpha}{2^{(k+1)/s}}  \rbrace  \big| d \alpha \\
& \leq \int_0^\infty \alpha^{p-1} \big(  \frac{\alpha}{2^{(k+1)/s}} \big)^{-p_2} \, K_2^{p_2} \big\| \big( \sum_{I \in \ii I_k} \frac{|\langle f, \phi_I^1 \rangle|^2}{|I|} \cdot \one_I  \big)^\frac{1}{2} \big\|_{L^{p_2}(L^Q)}^{p_2} \, d \alpha \\
&= 2^{p_2 \frac{k+1}{s}} K_2^{p_2} \int_0^\infty \alpha^{p-p_2-1} \int_{S(k, \alpha)} \| S_{\ii I_k} f(x)   \|_{L^Q}^{p_2} \, d x d \alpha. 
\end{align*}}

Changing the order of integration will not be helpful in this case because the collections on intervals $\ii I_0, \ii I_1, \ldots$ depend on the variable $\alpha$, and lower bounds for $S_{\I_k}$ independent on $\alpha$ are not available. Instead, we use the pointwise inequality $\ds \|S_{\ii I_k} f(x)\|_{L^Q} \leq \frac{\alpha}{C^{k-1}}$ from \eqref{eq:SF-est-S:k:aplh}. Recalling also the definition of $S(k, \alpha)$, we have 
{\fontsize{11}{11}\begin{align*}
&  2^{p_2 \frac{k+1}{s}} K_2^{p_2} \int_0^\infty \alpha^{p-p_2-1} \int_{S(k, \alpha)} \| S_{\ii I_k} f(x)   \|_{L^Q}^{p_2} \, d x d \alpha \\
&\leq 2^{p_2 \frac{k+1}{s}} K_2^{p_2} \int_0^\infty \alpha^{p-p_2-1}  |S(k, \alpha)| \big( \frac{\alpha}{C^{k-1}}  \big)^{p_2} d \alpha \\
& \leq 2^{p_2 \frac{k+1}{s}} K_2^{p_2} C^{-(k-1)p_2} \int_0^\infty \alpha^{p-1} \big| \lbrace x: \| S_{\ii I}f(x)  \|_{L^Q}> \frac{\alpha}{C^k}  \rbrace d 
\alpha \big|.
\end{align*}}

Making a change of variable we obtain
{\fontsize{11}{11}\begin{align*}
&  \int_{0}^\infty \alpha^{p-1}  \big| \lbrace x : \|  F_k(x)  \|_{L^Q}> \frac{\alpha}{2^{(k+1)/s}}  \rbrace  \big| d \alpha \\
& \leq 2^{p_2 \frac{k+1}{s}} K_2^{p_2} C^{-(k-1)p_2} C^{k p} \int_0^\infty \lambda^{p-1} \big| \lbrace x: \| S_{\ii I}f(x)  \|_{L^Q}> \lambda \rbrace \big| d 
\lambda\\
& \leq  \frac{1}{p} \,  2^{p_2 \frac{k+1}{s}} C^{p_2}  K_2^{p_2} C^{-k(p_2-p)}\|S_{\ii I} f\|_{L^p(L^Q)}^p.
\end{align*}}

Now it remains to put everything together and to sum in $k \geq 0$: due to \eqref{eq:interp:subadd},
\begin{align*}
\|\sum_{I \in \ii I} \langle f, \phi_I^1  \rangle \phi_I^2\|_{L^p(L^Q)}^p \leq \Big( \frac{p}{p-p_1} 2^\frac{p_1}{s} K_1^{p_1} + 2^{p_2/s} K_2^{p_2} C^{p_2} \sum_{k \geq 1} \big( 2^{p_2/s} C^{-(p_2-p)}  \big)^k \big)  \, \|S_{\ii I} f\|_{L^p(L^Q)}^p.
\end{align*}

Since $p_2-p>0$, if $C$ is large enough so that $\ds  2^{p_2/s} C^{-(p_2-p)} <1$ (which is equivalent to $C>2^{p_2/s(p_2-p)}$), the series above is finite. We obtain in this way \eqref{inter2} with 
\[
K^p \lesssim_{p_1, p_2, s} K_1^{p_1} + K_2^{p_2}. 
\]

\end{proof}
\end{proposition}

\begin{observation}
In the statement of Proposition \ref{prop:interpolation}, we could allow $K_1$ and $K_2$ to depend on the collection $\I$, which will yield an upper bound for $K$ that also depends on $\I$. Thus, assuming that
\begin{equation}\label{inter1}
\| (\sum_{I \in \ii I} \langle f, \phi^1_I\rangle \phi^2_I ) \one_E \|_{L^{p_j,\infty}(L^Q)} \leq K_j (\I) \,  \big\| \big( \sum_{I\in \ii I} \frac{ |\langle f, \phi^1_I\rangle|^2  } {|I|} \one_I \big)^{1/2} \big\|_{L^{p_j}(L^Q)}
\end{equation}
holds for all collections $\I$ of dyadic intervals, for $0<p_1<p<p_2<\infty$, we deduce the strong bound
\begin{equation}\label{inter2}
\| (\sum_{I \in \ii I} \langle f, \phi^1_I\rangle \phi^2_I ) \one_E\|_{L^p(L^Q)} \leq K(\I) \, \big\| \big( \sum_{I\in \ii I} \frac{ |\langle f, \phi^1_I\rangle|^2  } {|I|} \one_I \big)^{1/2} \big\|_{L^p(L^Q)},
\end{equation}
where this time $\ds K(\I) \lesssim \big( \sup_{\I' \subseteq \I } \, K_1 (\I')^{p_1} +  \sup_{\I' \subseteq \I } \, K_2(\I')^{p_2} \big)^\frac{1}{p} $.
\end{observation}

\section{Proof of Theorem \ref{mainth} in the general case}
\label{sec:proof-main-thm-gen-case}
Recall that our goal is to prove that
\begin{equation}\label{a1}
\| f \|_{L^P(L^Q)} \lesssim \| S f \|_{L^P(L^Q)}
\end{equation}
where the $d$-tuple $P = (p_1, ..., p_d)$ and the $n$-tuple $Q = (q_1, ..., q_n)$ satisfy $0 < P, Q < \infty$ componentwise. Recall also that the $N$-parameter square function $S $ is defined by
$$ S  := S_{d_1} \otimes ... \otimes S_{d_N}$$
while $d_1 + ... +d_N =d$. So far we have proved this in the particular situation when $d_1 = ... = d_N = 1$. The goal of this section is to explain that similar ideas
can handle the general case as well. 
First of all, let us observe that using a similar inductive argument to that in Section $2$, it is enough to prove the particular case when $N=1$. In other words, from now on, our square function $S f$ is 
a one parameter square function in $\R^d$ and the task is to prove multiple vector-valued, mixed norm estimates for it, in the form of
\begin{equation}\label{a2}
\| f \|_{L^P(L^Q)} \lesssim \| S_d f \|_{L^P(L^Q)}.
\end{equation}

It is now important to observe that when $p_1 = ... = p_d = p$, then (\ref{a2}) becomes a multiple vector-valued $L^p(\R^d)$ estimate, which can be proved exactly as in the one dimensional case $d=1$ treated before.
This is because all of our previous arguments have natural higher dimensional analogues. Instead of doing analysis with dyadic intervals, one does analysis with dyadic cubes of the corresponding dimension, in precisely the same way.

It will be more convenient to modify a bit the notation, in order to obtain a statement more suitable to the upcoming inductive argument. We will think of the Euclidean space $\R^d$ as being decomposed into
\begin{equation}\label{b1}
\R^d = \R^{n_1} \times ... \times \R^{n_M}
\end{equation}
and consequently the mixed norm space $L^P(\R^d)$ being unfolded as
\begin{equation}\label{b2}
L^P(\R^d) = L^{p_1}(\R^{n_1}) (L^{p_2}(\R^{n_2}) ( ... (L^{p_M}(\R^{n_M})) ... )).
\end{equation}
In other words, we implicitly assume that the first $n_1$ indices of the $d$-tuple $P$  are all equal to $p_1$, the next $n_2$ indices of $P$ are all equal to $p_2$, and so on, until the last $n_M$ set of indices of $P$ which are all equal to $p_M$.

The plan is to prove the corresponding (\ref{a2}) by induction with respect to the parameter $M$. As pointed out before (\ref{a2}) is already known when $M=1$ and we aim to show that it is also true for $M=d$, when all the entries
of $P$ are possibly different from each other.  

As in the one dimensional case, it is not difficult to see that things can be reduced to proving a discrete analogue of  (\ref{a2}) namely
\begin{equation}\label{a3}
\| \sum_{R\in \r} \langle f, \phi^1_R\rangle \phi^2_R \|_{L^P(L^Q)} \lesssim \Big\| \big( \sum_{R\in\r} \frac{ |\langle f, \phi^1_R\rangle|^2  } {|R|} \one_R \big)^{1/2} \Big\|_{L^P(L^Q)}.
\end{equation}
The families $(\phi^1_R)_R$ and $(\phi^2_R)_R$ in (\ref{a3}) are two lacunary families, $L^2$ normalized, indexed by a finite collection $\r$ of dyadic cubes in $\R^d$. And also as in the one dimensional case, the statement of Observation \ref{obs} remains
valid, in the sense that the two families of functions may depend on the implicit variables $(s_1, ..., s_n)$ of the space $L^Q$.

Using a higher dimensional analogue of (\ref{desco}) we decompose each $\phi^2_R$ as
\begin{equation}\label{desco1}
\phi^2_R = \sum_{k=0}^{\infty} 2^{ - \# k} \phi^2_{R, k} = : \sum_{k=0}^{\infty} 2^{-(\#/2) k} \widetilde{\phi}^2_{R, k} 
\end{equation}
where
$$\supp (\phi^2_{R,k}) \subseteq 2^k R$$
as before and where $\#$ is arbitrarily large. Using this in (\ref{a3}), it will be enough to show
\begin{equation}\label{a4}
\| \sum_{R\in \r} \langle f, \phi^1_R\rangle \widetilde{\phi}^2_{R, k} \, \|_{L^P(L^Q)} \lesssim 2^{ L k} \Big\| \big( \sum_{R\in\r} \frac{ |\langle f, \phi^1_R\rangle|^2  } {|R|} \one_R \big)^{1/2} \Big\|_{L^P(L^Q)}
\end{equation}
for some large but fixed number $L$. The main case is when $k=0$ and we will concentrate on it from now on (by this we mean that the general case follows by standard modifications as in the one dimensional situation). Then (\ref{a4}) reads as
\begin{equation}\label{a5}
\| \sum_{R\in \r} \langle f, \phi^1_R\rangle \widetilde{\phi}^2_{R, 0} \|_{L^P(L^Q)} \lesssim  \Big\| \big( \sum_{R\in\r} \frac{ |\langle f, \phi^1_R\rangle|^2  } {|R|} \one_R \big)^{1/2} \Big\|_{L^P(L^Q)}.
\end{equation}
We think of the dyadic cubes $R$ as being of the form
$$R = R_1 \times ... \times R_M$$
to match the decomposition (\ref{b1}), where each $R_j$ is a dyadic cube in $\R^{n_j}$ of the same side length as $R$ itself for $1\leq j\leq M$.

Following the same earlier strategy for the estimate (\ref{a5}), one needs in fact to prove a more localized variant of it given by
\begin{equation}\label{a6}
\| \sum_{R\in \r} \langle f, \phi^1_R\rangle \widetilde{\phi}^2_{R, 0} \one_E \|_{L^P(L^Q)} \lesssim  \Big\| \big( \sum_{R\in\r} \frac{ |\langle f, \phi^1_R\rangle|^2  } {|R|} \one_R \big)^{1/2} \Big\|_{L^P(L^Q)}  \cdot \big( \size_{\r_1} \one_E \big)^{1/p - \epsilon}
\end{equation}
where 
$$\r_1 := \{ R_1 : R=R_1 \times ... \times R_M \in \r \}$$
and  $\size_{\r_1} \one_E$ is the corresponding $n_1$-dimensional {\it size} generalizing naturally the one dimensional (\ref{sss}). In (\ref{a6}) the set $E$ is an arbitrary measurable subset of $\R^{n_1}$.

The plan is to prove (\ref{a6}) by induction with respect to the parameter $M$. Notice that
when $M=1$, then $\r_1 = \r$ and the corresponding (\ref{a6}) is known, as we pointed out before (its proof is identical to the one in the one dimensional case). In particular, all one has to do
is to prove that the case $M-1$ implies the case $M$, for every $M\geq 2$. We claim that this can de done by an argument similar to the one used earlier in the proof of ``$\P(n-1)$ implies $\P(n)$" (see Section \ref{sec:inductive-section}).

First of all, we like to see the left hand side of (\ref{a6}) as being
\begin{equation}\label{a7}
\| \sum_{R\in \r} \langle f, \phi^1_R\rangle \widetilde{\phi}^2_{R, 0} \one_E \|_{L^{p_1}(L^{\widetilde{P}}(L^Q))}
\end{equation}
where for $P = (p_1, ..., p_M)$ we define $\widetilde{P} := (p_2, ..., p_M)$. As before, by interpolation it would be enough to estimate the weaker analogue of it, namely
\begin{equation}\label{a8}
\| \sum_{R\in \r} \langle f, \phi^1_R\rangle \widetilde{\phi}^2_{R, 0} \one_E \|_{L^{p_1, \infty}(L^{\widetilde{P}}(L^Q))}
\end{equation}
by the same right hand side of (\ref{a6}). As explained previously, we dualize the $L^{p_1, \infty}$ quasi-norm through $L^s$, where $s$ is a positive real number smaller than all the entries of $P$, of $Q$, and also than $1$. By scale invariance (in the ambient space $\R^d$) this amounts
to prove that for every subset $F \subset \R^{n_1}$ with $|F| = 1$ there exists a major subset $\widetilde{F}\subseteq F$ with $|\widetilde{F}| \geq 1/2$ so that

\begin{equation}\label{a9}
\| \sum_{R\in \r} \langle f, \phi^1_R\rangle \widetilde{\phi}^2_{R, 0} \one_E  \one_{\widetilde{F}} \|_{L^s (L^{\widetilde{P}}(L^Q))} \lesssim R H S (\ref{a6}).
\end{equation}
The subset $\widetilde{F}$ is defined as usual by $\widetilde{F} := F \setminus \Omega$ for a certain exceptional set $\Omega\in \R^{n_1}$. This exceptional set is constructed as before with the only difference that the corresponding $\Omega_k$ are now given by
\begin{equation}\label{a10}
\Omega_k : = \{ x_1 \in \R^{n_1} : \|S f (x_1)\|_{L^{\widetilde{P}}(L^Q)} > C 2^{10 k /p_1} \| S f \|_{L^P(L^Q)} \}.
\end{equation}
In the above (\ref{a10}), by $S f$ one denotes the discrete square function given by the inner expression in the right hand side of (\ref{a3}). Also, we now think of a generic variable in $\R^d$ as being of the form
$(x_1, ..., x_M)$ with $x_j \in \R^{n_j}$ for $1\leq j\leq M$. In particular, $S f (x_1)$ can be thought of as a function depending on the rest of the variables $(x_2, ..., x_M)$ in an obvious way
$$S f (x_1) (x_2, ..., x_M) := S f (x_1, x_2, ..., x_M).$$
To estimate (\ref{a9}) one needs to perform (again) two carefully designed stopping times. The second one involves {\it averages} over dyadic cubes, and it is essentially a higher dimensional analogue of the one before. The first one on the other hand, selects maximal dyadic cubes $R_1^0$
in $\R^{n_1}$ for which the corresponding averages
\begin{equation}\label{a11}
\frac{1}{|R^0_1|^{1/p_1}}
\Big\| \big( \sum_{R\in\r : R_1\subseteq R_1^0} \frac{ |\langle f, \phi^1_R\rangle|^2  } {|R|} \one_R \big)^{1/2} \Big\|_{L^{p_1}(L^{\tilde P}(L^Q))}
\end{equation}
are large, also as in the one dimensional multiple vector-valued case. The way one uses these two together is similar to the way explained in the earlier ``$\P(n-1)$ implies $\P(n)$" situation. 
At some point, exactly as before, one uses H\"{o}lder locally, to be able to rely on the induction hypothesis (as in the previous (\ref{eq19})) in the particular case when $p_1 = p_2$.  More precisely, this amounts
to estimate expressions of type
$$\| \sum_{R} \langle f, \phi^1_R\rangle \widetilde{\phi}^2_{R, 0} \one_E \|_{L^{p_2}(L^{\widetilde{P}}(L^Q))}$$
locally, and here the induction hypothesis can be applied since the new $P$ tuple now is $P = (p_2, p_2, ..., p_M)$ and in particular, one can think of $\R^d$ as being split as $\R^d = \R^{n_1+n_2} \times ... \times \R^{n_M}$
and this contains now only $M-1$ factors. There are only two observations that one needs to make in order to realize that the earlier argument goes through smoothly in our case as well.

The first is that the John-Nirenberg inequality is still available in this context. More explicitly, this means that the supremum over $R^0_1$ of averages of type 
$$
\frac{1}{|R^0_1|^{1/p_2}}
\Big\| \big( \sum_{R\in\r : R_1\subseteq R_1^0} \frac{ |\langle f, \phi^1_R\rangle|^2  } {|R|} \one_R \big)^{1/2} \Big\|_{L^{p_2}(L^{\tilde P}(L^Q))},
$$
which appear naturally after one applies the induction, is controlled by the corresponding supremum of averages of type
$$
\frac{1}{|R^0_1|^{1/p_1}}
\Big\| \big( \sum_{R\in\r : R_1\subseteq R_1^0} \frac{ |\langle f, \phi^1_R\rangle|^2  } {|R|} \one_R \big)^{1/2} \Big\|_{L^{p_1}(L^{\tilde P}(L^Q))},
$$
which are the ones needed to capitalize on the stopping time procedure. To prove this, one just has to observe that the above inner expressions can also be seen as 
$$
 \sum_{R \in\r : R_1\subseteq R_1^0} \frac{ |\langle f, \phi^1_R\rangle|^2  } {|R|} \one_R (x_1, x_2, ..., x_M) = \sum_{R_1\subseteq R^0_1} \frac{|a_{R_1}(x_2, ..., x_M)|^2}{|R_1|} \one_{R_1} (x_1)
$$ 
where in general
$$a_{C} (x_2, ..., x_M) := \big( \sum_{R: R_1 = C} \frac{ |\langle f, \phi^1_R\rangle|^2  } {|R_2 \times \ldots \times R_M|} \one_{R_2 \times \dots\times R_M}(x_2, ..., x_M) \big)^\frac{1}{2}$$
and after that to realize that $BMO$ expressions of type
\begin{equation}\label{bmo}
\sup_{C_0} \frac{1}{|C_0|^{1/q}} \Big\| \big( \sum_{ C \subseteq C_0} \frac{|a_C|^2}{|C|} \one_C \big)^\frac{1}{2} \Big \|_{L^q(B)}
\end{equation}
are all equivalent to each other for every $0<q<\infty$ even when $B$ is a quasi-Banach lattice.

And the second observation is that 
$$\size_{\r_{1\times 2}} \one_E \lesssim \size_{\r_1} \one_E$$
as one can easily check. By $\r_{1\times 2}$ one means 
$$\r_{1\times 2}: = \{ R_1\times R_2 : R = (R_1, R_2, ..., R_M) \in \R \}$$
and they appear naturally after the application of the induction hypothesis in $\R^d = \R^{n_1+n_2} \times \ldots \times \R^{n_M}$. This concludes our proof of the weaker estimate (\ref{a9}).

After that the induction argument works exactly as before, allowing one to complete the proof of the desired discrete estimate (\ref{a6}).

\section{Connections to  weighted theory and extrapolation}
\label{sec:weights-extrap}

In the present section we discuss a certain weighted version of inequality \eqref{doua}, which eventually yields an alternative proof of Theorem \ref{mainth}, upon adapting existing extrapolation results. Assuming such a weighted estimate, in Section \ref{sec:extrapolation-proof}, we detail this proof by extrapolation. In the second part, Section \ref{sec:weighted-result-hel-method}, we review the weighted estimates (which are indispensable for extrapolation) and provide a proof for them based on a sparse domination result implied by the helicoidal method.

A weighted, scalar version of Theorem \ref{mainth} can be formulated in the following way:  if $f$ is a Schwartz function and $w$ is ``regular enough",
\begin{equation}
\label{eq:weighted-H-space}
\|f\|_{L^p(w)} \lesssim \| S f\|_{L^p(w)}.
\end{equation}
For $0< p \leq 1$, this inequality is related to the theory of weighted Hardy spaces and it was stated in \cite{DingHanLu-weightedHardySpaces}. There, the authors study the boundedness of singular integrals on such spaces, which was known previously under more stringent conditions on the weights (they were assumed to be $A_1$ weights). In \cite{DingHanLu-weightedHardySpaces}, a theory of weighted Hardy spaces and boundedness of singular integrals is developed for $A_\infty$ weights. Central to their theory is the inequality \eqref{eq:weighted-H-space}, which is stated for $A_\infty$ weights. Starting from this and using a certain type of extrapolation (regarding collections of pairs of functions, rather than operators, and $A_\infty$ weights), we recover the multiple vector-valued results of Theorem \ref{mainth}; the mixed-normed estimates are obtained through a generalization of a result of Kurtz \cite{KurtzExtrapMixedNormedSpaces}.

On the other hand, we will see once again that a local estimate similar to \eqref{discrete-2} and a change in the direction of the stopping time will yield a (multiple vector-valued) sparse estimate, and in consequence, also (multiple vector-valued) weighted estimates, in the one-parameter case. The weighted estimates obtained in this way are similar to \eqref{eq:weighted-H-space} and to those of \cite{DingHanLu-weightedHardySpaces}, and hence they are interconnected to weighted Hardy spaces.

Before proceeding, we briefly recall a few definitions and results about weights: if $1<p < \infty$, the measurable function $w: \rr R^m \to [0, \infty]$ belongs to the $A_p(\rr R^m)$ class provided
\[
[w]_{A_p}:= \sup_{\substack{Q \subset \rr R^m \\ Q \text{  cube}}}\big(\aver{Q} w(x) dx \big)  \, \big(\aver{Q} w^{1-p'}(x) dx \big)^{p-1} <+\infty.
\]
If $p=1$, then $w \in A_1(\rr R^m)$ provided there exists a constant $C$ such that $Mw(x) \leq C\, w(x)$ for almost every $x \in \rr R^m$. Then $A_\infty(\rr R^m)$ is defined as 
\[
A_\infty(\rr R^m):=\bigcup_{1 \leq p <\infty} A_p(\rr R^m).
\]
For the classes $A_{p, \, Rectangle}(\rr R^{d_1} \times \ldots \times \rr R^{d_N})$, the collection of cubes in replaced by the collection of rectangles with sides parallel to the coordinate axes, and in the case $p=1$, the Hardy-Littlewood maximal function is replaced by the \emph{strong maximal function} $M_S$. For $p>1$, it is well known that $w(x_1, \ldots, x_N) \in A_{p, \, Rectangle}(\rr R^{d_1} \times \ldots \times \rr R^{d_N})$ if and only if
\[
w(\cdot, x_2, \ldots, x_N) \in A_p(\rr R^{d_1}), \ldots, w( x_1, \ldots, x_{N-1}, \cdot) \in A_p(\rr R^{d_N}),
\]
uniformly with respect to the fixed variables.

\subsection{Weighted Hardy spaces and extrapolation}
\label{sec:extrapolation-proof}
Let $0<p<\infty$. If $w \in A_{\infty}(\rr R^m)$, then the \emph{weighted Hardy space} $H^p_w$ consists of
\begin{equation}
\label{def:weightedHardySpaces}
H^p_w := \lbrace f : \rr R^m \to \rr C : S_m(f) \in L^p_w(\rr R^m)  \rbrace.
\end{equation}

Setting $\ds \|f\|_{H_w^p} :=\|S_m(f)\|_{L^p_w}$, $H_w^p$ becomes a quasi-Banach space, for which we have, whenever $s \leq \min(p, 1)$
\[
\|f+g\|_{H_w^p}^s \leq \|f\|_{H_w^p}^s + \|g\|_{H_w^p}^s.
\]

By making use of a certain discrete Calder\'on reproduction formula, it was shown in \cite[Theorem~3.5]{DingHanLu-weightedHardySpaces} that, for any $w \in A_\infty(\rr R^m)$ and any $0<p \leq 1$, 
\begin{equation}
\label{eq:one-param-ineq}
\|f\|_{L^p_w(\rr R^m)} \leq C \| f \|_{H^p_w(\rr R^m)}=C\,\|S_m(f)\|_{L^p_w(\rr R^m)}.
\end{equation}

The method of the proof doesn't immediately generalize to the case $p>1$. Instead, in this situation the $L^p_w(\rr R^m)$ boundedness (which requires the stronger condition that $w \in A_p$) of the square function $S_m$ is invoked to deduce, by means of duality, an estimate similar to \eqref{eq:one-param-ineq}. Hence, for $p>1$, Ding et al. \cite{DingHanLu-weightedHardySpaces} state the inequality \eqref{eq:one-param-ineq} only for weights $w \in A_p$.

Alternatively, one can use the $A_\infty$ extrapolation developed in \cite{multi-extrap-CU-martell-perez} (similarly, see \cite[Corollary~3.15]{extrapolation-CUMP}) applied to the pairs of functions $(f, S_m (f))$. This will imply that \eqref{eq:one-param-ineq} is valid for any $0<p<\infty$, and for any $w \in A_\infty(\rr R^m)$. The same extrapolation result yields multiple vector-valued weighted inequalities: for any $0<p<\infty$, any $n$-tuple $Q$, and any weight $w \in A_\infty(\rr R^m)$, 
\begin{equation}
\label{eq:vv-Hardy-spaces}
\| f\|_{L^p(L^Q)(w)} \leq C \,\|S_m( f)\|_{L^p(L^Q)(w)}.
\end{equation}

Theorem 3.5 in \cite{DingHanLu-weightedHardySpaces} remains valid in the context of multi-parameter Hardy spaces, and Theorem 2.1 in \cite{multi-extrap-CU-martell-perez} holds for weights associated to Muckenhoupt bases. As a result, the multi-parameter multiple vector-valued inequality holds:
\begin{equation}
\label{eq:vv-Hardy-spaces-multiparam}
\| f\|_{L^p(L^Q)(w)} \leq C\|S_{d_1} \otimes\ldots \otimes  S_{d_N} ( f)\|_{L^p(L^Q)(w)},
\end{equation}
where $0<p<\infty$, $Q=(q_1, \ldots, q_n)$ with $0<q_j<\infty$ for all $1 \leq j \leq n$, and $w \in A_{\infty, Rectangle}(\rr R^{d_1} \times \ldots \times \rr R^{d_N})$.

In order to obtain the full mixed-norm estimates of Theorem \ref{mainth}, we need an extrapolation result from \cite{KurtzExtrapMixedNormedSpaces} suited for mixed-norm spaces. The result extends without any important modification to pairs of functions, in which case the operator $T$ is being disregarded. Once inequality \eqref{eq:vv-Hardy-spaces-multiparam} is deduced as above, the plan is to apply it to product weights and deduce the mixed norm estimates from Theorem \ref{thm-extrapolation-mixed-Kurtz} below.

We have the following reformulation of Kurtz's result, in a slightly more general setting, although the proof remains the same:

\begin{theorem}[Similar to Theorem 2 of \cite{KurtzExtrapMixedNormedSpaces}]
\label{thm-extrapolation-mixed-Kurtz}
Let $0<s_0<\infty$ and assume that there exists $s_0<s<\infty$ such that
\begin{equation}
\label{eq:hyp-Kurtz-extrap-mixed-rescaled}
\int_{\rr R^{d_1} \times \rr R^{d_2}} |f(x, y)|^s w(x, y) \,dy \,dx \leq C \int_{\rr R^{d_1} \times \rr R^{d_2}} |g(x, y)|^s w(x, y) \,dy \,dx
\end{equation}
for all pairs $(f, g)$ in a certain collection of functions $\ic F$, and for all $w \in A_{\frac{s}{s_0}, \, Rectangle}(\rr R^{d_1} \times \rr R^{d_2})$, with a constant depending only on $[w]_{A_{\frac{s}{s_0}, \, Rectangle}}$. Then for any $s_0 <p, q <\infty$, and any weights $w(x, y)$ of the type $w(x,y)=u(x) \, v(y)$ such that 
\[
u^\frac{p}{q} \in A_{\frac{p}{s_0}}(\rr R^{d_1}), \qquad v \in A_{\frac{q}{s_0}}(\rr R^{d_2}),
\]
we have
\[
\int_{\rr R^{d_1}} \big( \int_{ \rr R^{d_2}} |f(x, y)|^q w(x, y) \,dy \big)^\frac{p}{q} \,dx \leq C( [u]_{A_{\frac{p}{s_0}}}, [v]_{A_{\frac{q}{s_0}}} ) \int_{\rr R^{d_1}} \big( \int_{ \rr R^{d_2}}  |g(x, y)|^q w(x, y) \,dy \big)^\frac{p}{q} \,dx.
\]

In particular, if $w(x,y)\equiv 1$, mixed-norm estimates are implied by extrapolation, once the weighted result \eqref{eq:hyp-Kurtz-extrap-mixed-rescaled} is known.
\end{theorem}

\begin{remark}
In \cite{KurtzExtrapMixedNormedSpaces}, one is in fact looking for a necessary and sufficient conditions on weights $w(x,y)$ so that the strong maximal function $M_S$ satisfies 
\[
\int_{\rr R^{d_1}} \big( \int_{ \rr R^{d_2}} |M_S f(x, y)|^q w(x, y) \,dy \big)^\frac{p}{q} \,dx \leq C \int_{\rr R^{d_1}} \big( \int_{ \rr R^{d_2}}  |f(x, y)|^q w(x, y) \,dy \big)^\frac{p}{q} \,dx.
\]
While a necessary condition was found (the classes $A_p(A_q)$ from \cite[Definition~2]{KurtzExtrapMixedNormedSpaces}), sufficiency is proved only in the particular case of product weights $w(x, y)=u(x) \, v(y)$. Since we are mainly interested in the unweighted, multiple vector-valued case, we do not elaborate on the properties of the class of weights $A_p(A_q)$, but instead focus on the extrapolation result, which is also known to be true only for product weights.

We also don't keep track how the constants appearing in the inequalities above depend on the weights involved or on their characteristics.
\end{remark}

Next, we generalize Theorem \ref{thm-extrapolation-mixed-Kurtz} to mixed-norm $L^p$ spaces involving $\kappa$ variables (with $\kappa \geq 2$), and $A_\infty$ weights.

\begin{theorem}
\label{thm-extrap-kurtz-A-infty}
Assume there exists some $0<s<\infty$ so that
{\fontsize{9}{9}\begin{equation} \label{eq:hyp-thm-extrap-kurtz-A-infty}
\int_{\rr R^{d_1} \times \ldots \times \rr R^{d_\kappa}} |f(x_1, \ldots, x_\kappa)|^s w(x_1, \ldots, x_\kappa) d x_1 \ldots d x_\kappa \leq C \, \int_{\rr R^{d_1} \times \ldots \times \rr R^{d_\kappa}} |g(x_1, \ldots, x_\kappa)|^s w(x_1, \ldots, x_\kappa) d x_1 \ldots d x_\kappa
\end{equation}}for all $w \in A_{\infty, Rectangle}(\rr R^{d_1} \times \ldots \times \rr R^{d_\kappa})$ and for all pairs of functions $(f, g) \in \ic F$. Then for any $0<p_1, \ldots, p_\kappa <\infty$ and for any weight $w(x_1, \ldots, x_\kappa)=w_1(x_1) \cdot \ldots \cdot w_\kappa(x_\kappa)$ so that $\ds w_l^{\frac{p_l}{p_\kappa}} \in A_\infty(\rr R^{d_l})$ for all $1 \leq l \leq \kappa$, we have 
{\fontsize{9}{9}
\begin{align}
\label{eq:ineq-A-infty-extrap-kurtz}
&\big(\int_{\rr R^{d_1}} \ldots \big( \int_{ \rr R^{d_\kappa}} |f(x_1, \ldots, x_\kappa)|^{p_\kappa} w(x_1, \ldots, x_\kappa) d x_\kappa\big)^\frac{p_{\kappa-1}}{p_\kappa} \ldots d x_1 \big)^\frac{1}{p_1} \\
& \qquad \leq C \, \big(\int_{\rr R^{d_1}} \ldots \big( \int_{ \rr R^{d_\kappa}} |g(x_1, \ldots, x_\kappa)|^{p_\kappa} w(x_1, \ldots, x_\kappa) d x_\kappa \big)^\frac{p_{\kappa-1}}{p_\kappa} \ldots d x_1 \big)^\frac{1}{p_1}, \nonumber
\end{align}
}
for all $(f, g) \in \ic F$.
\begin{proof}
We present a proof by induction over $\kappa$. If $\kappa=2$, the statement is a reformulation of Theorem \ref{thm-extrapolation-mixed-Kurtz}:
the assumption that $w_1^\frac{p_1}{p_2}, w_2 \in A_\infty$ will be rewritten so that $\ds w_1^\frac{p_1}{p_2} \in A_{\frac{p_1}{s_0}}, w_2 \in A_\frac{p_2}{s_0}$, for a suitable $0<s_0<\infty$.

Since $0<p_1, p_2 <\infty$ and $w_1^\frac{p_1}{p_2} \in A_\infty, w_2 \in A_\infty$, there exists $1 \leq s_1, s_2 < \infty$ so that 
\[
w_1^\frac{p_1}{p_2} \in A_{s_1}, \qquad w_2 \in A_{s_2}.
\]
We pick $s_0$ with $0<s_0 \leq s$ with $\ds s_1 \leq \frac{p_1}{s_0}, \, s_2 \leq \frac{p_2}{s_0}$ (these conditions reduce to $\ds s_0 \leq \min ( \frac{p_1}{s_1}, \frac{p_2}{s_2}, s)$). Because the weight classes are nested, we have in this situation $\ds w_1^\frac{p_1}{p_2} \in A_{\frac{p_1}{s_0}}, w_2 \in A_\frac{p_2}{s_0}$.

The hypothesis \eqref{eq:hyp-thm-extrap-kurtz-A-infty} holds for all weights $w \in A_{\infty, Rectangle}$, and in particular also for $w \in A_{\frac{s}{s_0}, Rectangle}$; the inequality in \eqref{eq:ineq-A-infty-extrap-kurtz} then follows from Theorem \ref{thm-extrapolation-mixed-Kurtz}.
 
Next, we assume that the result holds true when $\kappa-1$ variables are involved and will prove it for $\kappa$ variables as well. We fix a $\kappa$-tuple $(p_1, \ldots, p_\kappa)$ and weights $w_1, \ldots, w_\kappa$ satisfying $\ds w_l^\frac{p_l}{p_\kappa}$ for all $1 \leq l \leq \kappa$. Denote by 
\[
F(x_1, x_2):= \big \|  f(x_1, \ldots, x_\kappa)  \big\|_{L^{p_3}_{x_3} \ldots L^{p_\kappa}_{x_{\kappa}}(w_3 \cdot \ldots \cdot w_\kappa)}, \quad G(x_1, x_2):= \big \|  g(x_1, \ldots, x_\kappa)  \big\|_{L^{p_3}_{x_3} \ldots L^{p_\kappa}_{x_{\kappa}}(w_3 \cdot \ldots \cdot w_\kappa)}.
\]

We want to show that 
{\fontsize{9}{9}\[
\int_{\rr R^{d_1}} \big( \int_{ \rr R^{d_2}} |F(x_1, x_2)|^{p_2} w_1^\frac{p_2}{p_{\kappa}}(x_1) w_2^\frac{p_2}{p_{\kappa}}(x_2) \,dx_2 \big)^\frac{p_1}{p_2} \,dx_1 \leq C \int_{\rr R^{d_1}} \big( \int_{ \rr R^{d_2}}  |G(x_1, x_2)|^{p_2} w_1^\frac{p_2}{p_{\kappa}}(x_1) w_2^\frac{p_2}{p_{\kappa}}(x_2) \,dx_2 \big)^\frac{p_1}{p_2} \,dx_1, 
\]}
given that $w_1^\frac{p_1}{p_{\kappa}} \in A_{\infty}$ and $w_2^\frac{p_2}{p_{\kappa}} \in A_{\infty}$.

If we denote 
\[
W(x_1, x_2):=w_1^\frac{p_2}{p_{\kappa}}(x_1) \, w_2^\frac{p_2}{p_{\kappa}}(x_2):=U(x_1) \, V(x_2), 
\]
we have $U^\frac{p_1}{p_2} \in A_\infty$ and $V \in A_{\infty}$. The problem is reduced to the case $\kappa=2$, and it remains to check that the hypothesis \eqref{eq:hyp-thm-extrap-kurtz-A-infty} is satisfied. That is, we need to check that there exists $0<s <\infty$ such that 
\begin{equation}
\label{eq:almost-there}
\int_{\rr R^{d_1} \times \rr R^{d_2}} |F(x_1, x_2)|^s \, W_0(x_1, x_2) d x_1d x_2 \leq C \, \int_{\rr R^{d_1} \times \rr R^{d_2}} |G(x_1, x_2)|^s \, W_0(x_1, x_2) d x_1 d x_2
\end{equation}
for all weights $W_0 \in A_{\infty, Rectangle} (\rr R^{d_1} \times \rr R^{d_2})$.

The case of $(\kappa -1)$ iterated Lebesgue spaces, applied to the tuple $(\tilde p_2, p_3, \ldots, p_\kappa)$ for some $0< \tilde p_2 <\infty$ yields, for weights of the form $w(x_1, x_2, \ldots, x_\kappa)=\tilde w_2(x_1, x_2) \cdot w_3(x_3) \cdot \ldots  \cdot w_\kappa (x_\kappa)$ so that $\ds w_l^{\frac{p_l}{p_\kappa}} \in A_\infty(\rr R^{d_l})$ for all $3 \leq l \leq \kappa$ and $\ds \tilde w_2^\frac{\tilde p_2}{p_\kappa}(x_1, x_2) \in A_\infty(\rr R^{d_1+d_2})$ the estimate
{\fontsize{9}{9}
\begin{align*}
&\big(\int_{\rr R^{d_1+d_2}} \ldots \big( \int_{ \rr R^{d_\kappa}} |f(x_1,x_2,  \ldots, x_\kappa)|^{p_\kappa} w(x_1, \ldots, x_\kappa) d x_\kappa\big)^\frac{p_{\kappa-1}}{p_\kappa} \ldots d x_1 d x_2 \big)^\frac{1}{\tilde p_2} \\
& \qquad \leq C \, \big(\int_{\rr R^{d_1+d_2}} \ldots \big( \int_{ \rr R^{d_\kappa}} |g(x_1, x_2, \ldots, x_\kappa)|^{p_\kappa} w(x_1, \ldots, x_\kappa) d x_\kappa \big)^\frac{p_{\kappa-1}}{p_\kappa} \ldots d x_1 \, d x_2 \big)^\frac{1}{\tilde p_2}.
\end{align*}
}

If the functions $f$ and $g$ and the weights $w_3, \ldots, w_\kappa$ are precisely those we started with, we obtain, for any weight $\tilde w_2$ so that $ \tilde w_2(x_1, x_2)^\frac{\tilde p_2}{p_\kappa} \in A_\infty(\rr R^{d_1+d_2})$ the estimate
\begin{equation}
\label{eq:almost-there-non-rectangle}
\int_{\rr R^{d_1} \times \rr R^{d_2}} |F(x_1, x_2)|^{\tilde p_2} \tilde w_2(x_1, x_2) d x_1d x_2 \leq C \, \int_{\rr R^{d_1} \times \rr R^{d_2}} |G(x_1, x_2)|^{\tilde p_2} \tilde w_2(x_1, x_2) d x_1 d x_2.
\end{equation}

We want \eqref{eq:almost-there} for some $0<s <\infty$ and all weights $W_0(x_1, x_2) \in A_{\infty, Rectangle} (\rr R^{d_1} \times \rr R^{d_2})$. Instead, the $(k-1)$ induction case yields the similar estimate \eqref{eq:almost-there-non-rectangle} for any $0 <\tilde p_2<\infty$ and any weight $\tilde w_2$ so that $\tilde w_2(x_1, x_2)^\frac{\tilde p_2}{p_\kappa} \in A_\infty(\rr R^{d_1+d_2})$. We get the desired estimate by choosing $s=p_\kappa$ and by noting that the class of weights for which supremum over rectangles is finite is a subcollection the class of weights for which supremum over cubes is finite:
\[
A_{\infty, Rectangle} (\rr R^{d_1} \times \rr R^{d_2}) \subset A_\infty(\rr R^{d_1+d_2}).
\]

\end{proof}
\end{theorem}

\subsubsection*{Proof of the main Theorem \ref{mainth}}
Now we want to deduce the general inequality
\[
\|f\|_{L^P(L^Q)} \lesssim \| S(f) \|_{L^P(L^Q)}.
\]

By extrapolating the scalar result of \cite{DingHanLu-weightedHardySpaces}, we obtain the multiple vector-valued estimate of \eqref{eq:vv-Hardy-spaces-multiparam}. Then we apply Theorem \ref{thm-extrap-kurtz-A-infty} in the case of $d=d_1+ \ldots + d_N$ variables, to obtain the mixed-norm, multiple vector-valued result.

\subsection{Obtaining the weighted result by using the helicoidal method}
\label{sec:weighted-result-hel-method}

As previously mentioned, we can obtain the weighted result directly from a sparse domination estimate, which follows from a \emph{local maximal inequality}. A similar strategy was used in \cite{BM3}.

\subsubsection{The Localization Lemma}
\label{sec:localization-lemma}
For the weighted result, it is more suitable to work with locally integrable functions than with characteristic functions, the reason being that the characteristic function cannot play the role of an $A_\infty$ weight.

We recall a few notations, for convenience:
\begin{notation}
If $\ic I$ is a collection of cubes in $\rr R^d$ and $I_0 \subseteq \rr R^d$ is a fixed dyadic cube, then 
\[ 
\ic I (I_0) := \lbrace I \in \ic I: I \subseteq I_0   \rbrace \qquad \text{and} \qquad \ic I ^+(I_0) := \ic I(I_0) \cup \lbrace I_0 \rbrace.
\]
For any cube $I \subset \rr R^d$, $\ci_I(x)$ denotes a function that decays fast away from $I$:
\begin{equation}
\label{def:ci}
\ci_I(x) := \big( 1+ \frac{\dist(x, I)}{|I|} \big)^{-M},
\end{equation}
where $M$ can be as large as we wish.
\end{notation}

\begin{remark}
For statements involving a weight $w \in A_{\infty}$, the decaying factor $M$ in the definition \eqref{def:ci} might depend on $w$. More exactly, if $w \in A_\infty$, then we know that $w \in A_{q_w}$ for some $q_w >1$; we will need, in certain situations, to make sure that $d\, q_w < M$.
\end{remark}

\begin{lemma}
\label{lemma:localization-lemma}
Let $0<p \leq 1$. Let $\mathcal I$ be a finite collection of dyadic squares in $\rr R^d$, $I_0$ a fixed dyadic square, $f:\rr R^d \to \rr C$ a Schwartz function, and $w$ a locally integrable, positive function. Then 
{\fontsize{10}{10}\begin{equation}
\label{eq:local-eq}
\big\| \sum_{I \in \mathcal{I}(I_0)} \langle f, \phi_I^1 \rangle \, \phi_I^2   \big\|_{L^p(w)}^p \lesssim \big( \sup_{J_1 \in \ic{I}(I_0)} \frac{1}{|J_1|^\frac{1}{p_1}}   \Big\| \big(  \sum_{\substack{I \in \ic{I}(I_0) \\ I \subseteq J_1}} \frac{|\langle f, \phi_I^1  \rangle|^2}{|I|} \cdot \one_I  \big)^\frac{1}{2}   \Big\|_{p_1}  \big)^p \, \big( \sup_{J_2 \in \ic{I}^+(I_0)} \frac{1}{|J_2|} \int_{\rr R} w \cdot \ci_{J_2} dx  \big) \cdot |I_0|,
\end{equation}}with an implicit constant independent on the collection $\ic I$ and on the functions $f$ and $w$. 
\begin{proof}
If $0<p <1$, then $\|  \cdot \|_p^p$ is subadditive. In this case, we have for some $0<\tau<\infty$
\[
\frac{1}{p}=1+\frac{1}{\tau}.
\]

First, we note that $\ds \big\| \sum_{I \in \mathcal{I}(I_0)} \langle f, \phi_I^1 \rangle \, \phi_I^2   \big\|_{L^p(w)} =\big\| \big( \sum_{I \in \mathcal{I}(I_0)} \langle f, \phi_I^1 \rangle \, \phi_I^2 \big) \cdot w^\frac{1}{p}   \big\|_{L^p}$. We let $v_1:=w$ and $v_2:=w^\frac{1}{\tau}$, so that 
\[
w^\frac{1}{p}=v_1 \cdot v_2 \qquad \text{and} \qquad \big\| \sum_{I \in \mathcal{I}(I_0)} \langle f, \phi_I^1 \rangle \, \phi_I^2   \big\|_{L^p(w)}= \big\| \big( \sum_{I \in \mathcal{I}(I_0)} \langle f, \phi_I^1 \rangle \, \phi_I^2  \big) \, v_1 \cdot v_2 \big\|_{L^p}.
\]

We also use the previous decomposition \eqref{desco} $\ds \phi_I^2(x):=\sum_{\ell \geq 0} 2^{-\frac{\ell \, M}{2}} \tilde{\phi}^2_{I, \ell}(x)$, so that it suffices to show instead of \eqref{eq:local-eq} the similar inequality, for every $\ell \geq 0$:
{\fontsize{9}{9}\begin{equation}
\label{eq:local-eq-ell}
\| ( \sum_{I \in \mathcal{I}(I_0)} \langle f, \phi_I^1 \rangle \, \tilde{\phi}_{I, \ell}^2 ) \, v_1 \cdot v_2  \|_{L^p}^p \lesssim 2^{10 \ell \, d \, p} \big( \sup_{J_1 \in \ic{I}(I_0)} \frac{1}{|J_1|^\frac{1}{p_1}}   \Big\| \big(  \sum_{\substack{I \in \ic{I}(I_0) \\ I \subseteq J_1}} \frac{|\langle f, \phi_I^1  \rangle|^2}{|I|} \cdot \one_I  \big)^\frac{1}{2}   \Big\|_{p_1}  \big)^p \, \big( \sup_{J_2 \in \ic{I}^+(I_0)} \frac{1}{| J_2|} \int_{\rr R^d} w \cdot \ci_{J_2} dx  \big)^p \cdot |I_0|.
\end{equation}}

We recall that the families $\big( \tilde{\phi}_{I, \ell}^2  \big)_{I \in \ic{I}}$ are all lacunary, $L^2$-normalized, and $\supp \, \tilde{\phi}_{I, \ell}^2 \subseteq 2^\ell I$ for all $I \in \ic{I}$. As before, we only present the case $\ell=0$, since the general case follows from almost identical arguments.

By H\"older's inequality and the fact that all the functions $\tilde{\phi}_{I, 0}^2$ are supported on $I \subseteq I_0$, we have
\begin{equation}
\label{eq:holder-p<1}
\big\| \big( \sum_{I \in \mathcal{I}(I_0)} \langle f, \phi_I^1 \rangle \, \tilde{\phi}_{I, 0}^2 \big) \, v_1 \cdot v_2  \big\|_{L^p} \lesssim \big\| \big( \sum_{I \in \mathcal{I}(I_0)} \langle f, \phi_I^1 \rangle \, \tilde{\phi}_{I, 0}^2 \big) \, v_1  \big\|_{L^1} \cdot \| v_2 \cdot \one_{I_0}\|_\tau.
\end{equation}

The first expression can be rewritten as
\[
\int_{\rr R^d} \big( \sum_{I \in \mathcal{I}(I_0)} \langle f, \phi_I^1 \rangle \, \tilde{\phi}_{I, 0}^2(x) \big) \, v_1(x) \cdot \overline{ g(x)} dx=  \sum_{I \in \mathcal{I}(I_0)} \langle f, \phi_I^1 \rangle \, \overline{ \langle  \overline {v_1} \cdot g, \tilde{\phi}_{I, 0}^2 \rangle},
\]
for a certain function $g \in L^\infty$ satisfying $\|g\|_\infty=1$. Next, we will introduce square functions in order to make use of John-Nirenberg inequality(\cite[Theorem2.7]{cw}):
{\fontsize{10}{10}\begin{align*}
& \vert   \sum_{I \in \mathcal{I}(I_0)} \langle f, \phi_I^1 \rangle \, \langle v_1 \cdot g, \tilde{\phi}_{I, 0}^2 \rangle  \vert =\Big\vert \int_{\rr R^d}  \sum_{I \in \mathcal{I}(I_0)} \frac{ \langle f, \phi_I^1 \rangle}{|I|^{1/2}} \cdot \one_I(x) \, \frac{ \langle v_1 \cdot g, \tilde{\phi}_{I, 0}^2 \rangle}{|I|^{1/2}} \cdot \one_I(x) dx   \Big\vert \\
& \lesssim \frac{1}{|I_0|^{1/2}}  \Big\| \big(  \sum_{\substack{I \in \ic{I}(I_0) \\ I \subseteq I_0}} \frac{|\langle f, \phi_I^1  \rangle|^2}{|I|} \cdot \one_I  \big)^\frac{1}{2}   \Big\|_{2} \cdot \frac{1}{|I_0|^{1/2}} \Big\| \big(  \sum_{\substack{I \in \ic{I}(I_0) \\ I \subseteq I_0}} \frac{| \langle v_1 \cdot g, \tilde{\phi}_{I, 0}^2 \rangle|^2}{|I|} \cdot \one_I  \big)^\frac{1}{2}   \Big\|_{2} \cdot |I_0| \\
&\lesssim \big( \sup_{J_1 \in \ic{I}(I_0)} \frac{1}{|J_1|^\frac{1}{p_1}}\Big\| \big(  \sum_{\substack{I \in \ic{I}(I_0) \\ I \subseteq J_1}} \frac{|\langle f, \phi_I^1  \rangle|^2}{|I|} \cdot \one_I  \big)^\frac{1}{2}   \Big\|_{p_1} \big) \cdot \big( \sup_{J_2 \in \ic{I}^+(I_0)} \frac{1}{|J_2|^\frac{1}{p_2}}  \Big\| \big(  \sum_{\substack{I \in \ic{I}(I_0) \\ I \subseteq J_2}} \frac{| \langle v_1 \cdot g, \tilde{\phi}_{I, 0}^2 \rangle|^2}{|I|} \cdot \one_I  \big)^\frac{1}{2}   \Big\|_{p_2, \infty} \big) \, |I_0|,
\end{align*}}for any $0<p_1, p_2 <\infty$.

Setting $p_2=1$ and using the $L^1 \mapsto L^{1, \infty}$ boundedness of the square function (see also \cite[Lemma 2.13]{cw}) we obtain
{\fontsize{10}{10}\[
\big\| \big( \sum_{I \in \mathcal{I}(I_0)} \langle f, \phi_I^1 \rangle \, \tilde{\phi}_{I, \ell}^2 \big) \, v_1  \big\|_{L^1} \lesssim \big( \sup_{J_1 \in \ic{I}(I_0)} \frac{1}{|J_1|^\frac{1}{p_1}}\Big\| \big(  \sum_{\substack{I \in \ic{I}(I_0) \\ I \subseteq J_1}} \frac{|\langle f, \phi_I^1  \rangle|^2}{|I|} \cdot \one_I  \big)^\frac{1}{2}   \Big\|_{p_1} \big)  \cdot \big( \sup_{J_2 \in \ic{I}^+(I_0)} \frac{1}{|J_2|} \int_{\rr R^d} v_1 \cdot \ci_{J_2} dx  \big) \cdot |I_0|.
\]}
Recalling that $v_1=w$ and $\ds \| v_2 \cdot \one_{I_0}  \|_\tau = \big( \frac{1}{|I_0|} \| w \cdot \one_{I_0}\|_1  \big)^\frac{1}{\tau} \cdot |I_0|^\frac{1}{\tau}$, the above estimate and \eqref{eq:holder-p<1} imply that 
\begin{align*}
&\big\| \big( \sum_{I \in \mathcal{I}(I_0)} \langle f, \phi_I^1 \rangle \, \tilde{\phi}_{I, 0}^2 \big) \, v_1 \cdot v_2  \big\|_{L^p}   \lesssim  
\big( \sup_{J_1 \in \ic{I}(I_0)} \frac{1}{|J_1|^\frac{1}{p_1}}\Big\| \big(  \sum_{\substack{I \in \ic{I}(I_0) \\ I \subseteq J_1}} \frac{|\langle f, \phi_I^1  \rangle|^2}{|I|} \cdot \one_I  \big)^\frac{1}{2}   \Big\|_{p_1} \big) \\ 
&  \qquad \cdot \big( \sup_{J_2 \in \ic{I}^+(I_0)} \frac{1}{|J_2|} \int_{\rr R^d} w \cdot \ci_{J_2} dx  \big) \cdot |I_0| \cdot \big( \frac{1}{|I_0|} \| w \cdot \one_{I_0}\|_1  \big)^\frac{1}{\tau} \cdot |I_0|^\frac{1}{\tau}.
\end{align*}
Raising the inequality to power $p$ we obtain exactly the inequality \eqref{eq:local-eq-ell} in the case $\ell=0$.

For $\ell \geq 1$, the difference will consist in replacing \eqref{eq:holder-p<1} by
\[
\big\| \big( \sum_{I \in \mathcal{I}(I_0)} \langle f, \phi_I^1 \rangle \, \tilde{\phi}_{I, \ell}^2 \big) \, v_1 \cdot v_2  \big\|_{L^p} \lesssim \big\| \big( \sum_{I \in \mathcal{I}(I_0)} \langle f, \phi_I^1 \rangle \, \tilde{\phi}_{I, \ell}^2 \big) \, v_1  \big\|_{L^1} \cdot 2^{\frac{\ell \, M_1}{\tau }}  \big( \frac{1}{|I_0|} \int_{\rr R^d} v_2^\tau \cdot \ci_{I_0} dx   \big)^{\frac{1}{\tau}}
\]
and using the $L^1 \mapsto L^{1, \infty}$ boundedness of the modified square function
\[
g \mapsto \big(  \sum_{\substack{I \in \ic{I}(I_0) \\ I \subseteq J_2}} \frac{| \langle g, \tilde{\phi}_{I, \ell}^2 \rangle|^2}{|I|} \cdot \one_I  \big)^\frac{1}{2},
\]
which satisfies the same $L^p$ estimates as the classical disretized square function of \cite{cw}, uniformly in $\ell \geq 0$.

The inequality stays true if $p=1$; in that case $\tau=\infty$ and there will not be a second term on the right hand side of \eqref{eq:holder-p<1}.
\end{proof}
\end{lemma}

\begin{remark}
This should be compared to the maximal inequality in Theorem 19 of \cite{BM3}.
\end{remark}

\subsubsection{The stopping time}
Further, Theorem 12 in \cite{BM3} explains how to deduce sparse estimates from a local estimate such as \eqref{eq:local-eq} of Lemma \ref{lemma:localization-lemma}. The procedure in \cite{BM3} is stated for averages of functions, but the same is true when averages of square functions are concerned. 

A similar algorithm, based on the helicoidal method, was used in \cite{BB}, to deduce a sparse domination by averages of localized square functions result.

\begin{theorem}
\label{thm-sparse-dom}
Let $\ic{I}$ be a collection of dyadic squares, $0<p<\infty$ and $w$ a positive, locally integrable function. Then, for any $\epsilon_p >0$ and any Schwartz function $f$, there exists a sparse collection $\ic S$ of cubes (which depends on the functions $f, w$, the exponent $p$) so that 
{\fontsize{10}{10}\begin{equation}
\label{eq:sparse-SF}
\| \big( \sum_{I \in \mathcal{I}} \langle f, \phi_I^1  \rangle \,  \phi_I^2\big) \cdot w^\frac{1}{p} \|_p^p \lesssim \sum_{Q \in S} \big(  \frac{1}{|Q|^\frac{1}{p_1}}\Big\| \big(  \sum_{\substack{I \in \ic{I} \\ I \subseteq Q}} \frac{|\langle f, \phi_I^1  \rangle|^2}{|I|} \cdot \one_I  \big)^\frac{1}{2}   \Big\|_{p_1} \big)^p  \cdot \big( \frac{1}{|Q|} \int_{\rr R^d} w^{1+\epsilon_p} \, \ci_Q  \, dx \big)^\frac{1}{1+\epsilon_p} \cdot |Q|.
\end{equation}}
However, if $0<p \leq 1$, the above inequality is true for $\epsilon_p=0$.

\begin{proof}
We briefly sketch the proof for completeness, first in the case $0<p \leq 1$: per usual, the collection $ \ic S:=\bigcup\limits_{k \geq 0} \ic S_k$, where the cubes in the sub-collection $\ic S_{k+1}$ are to be understood as the ``descendants" of the dyadic cubes in the previous generation $\ic S_k$:
\[
\ic S_{k+1}:=\bigcup_{Q \in \ic S_k} ch_{\ic S}(Q).
\]
 To every $Q \in \ic S$, we also associate a subcollection $\ic I_Q \subseteq \ic I$ of cubes so that 
\[
\ic I := \bigcup_{Q \in \ic S} \ic I_Q
\]
represents a partition of the initial collection $\ic I$.

The bottom-most collection $\ic S_0$ will consist of the maximal dyadic cubes of the collection $\ic I$:
\[
\ic S_0:= \lbrace Q \in \ic I: Q \text{  maximal with respect to inclusion }  \rbrace.
\]

Next, we assume that $\ic S_0, \ic S_1 $ up to $\ic S_k$ are known and we will show how to construct $\ic S_{k+1}$, and for every $Q_0 \in \ic S_k$, the collections $\ic I_{Q_0}$.

If $Q_0 \in \ic S_k$, then we define
\begin{align}
\label{def-excp-set-sparse}
E_{Q_0}:&=  \big \lbrace x  \in Q_0 : \Big( \sum_{\substack{ I \in \ic I \\ I \subseteq Q_0}} \frac{|\langle f, \phi_I^1 \rangle|^2}{|I|} \one_I(x) \Big)^\frac{1}{2} > C \, \frac{1}{|Q_0|^{\frac{1}{p_1}}} \Big\| \big( \sum_{\substack{ I \in \ic I \\ I \subseteq Q_0}} \frac{|\langle f, \phi_I^1 \rangle|^2}{|I|} \one_I \big)^\frac{1}{2}  \Big\|_{p_1}  \big \rbrace \\
& \cup \big\lbrace  x \in Q_0 :   M( w \cdot \ci_{I_0})(x) > C \, \frac{1}{|Q_0|} \int_{\rr R^d} w \cdot \ci_{Q_0}(y) dy  \big\rbrace \nonumber.
\end{align}

It is not difficult to see that, if we choose $C>0$ large enough, $\ds |E_{Q_0}| >\frac{|Q_0|}{2}$. Then $ch_{\ic S}(Q_0)$ will consist of a maximal covering of $E_{Q_0}$ by dyadic cubes:
\[
ch_{\ic S}(Q_0):=\lbrace Q \text{ dyadic cube  }: Q \subseteq E_{Q_0}, \text{ maximal with respect to inclusion}   \rbrace
\]
and also, as already stated, $\ds \ic S_{k+1}:= \bigcup_{Q_0 \in \ic S_k} ch_{\ic S}(Q_0)$.

On the other hand, for every $Q_0 \in \ic S_k$, we define
\[
\ic I_{Q_0}:= \lbrace I \in \ic I : I \subseteq Q_0, I \nsubseteq E_{Q_0}   \rbrace.
\]
In consequence, every $I \in \ic I_{Q_0}$ has the property that either it is disjoint from the intervals in $ch_{\ic S}(Q_0)$, or, if $Q \in ch_{\ic S}(Q_0)$ and $I \cap Q \neq \emptyset$, then necessarily $Q \subsetneq I$. This implies in particular that the localized square function 
\begin{equation}
\label{def:localized-square-function}
S_{\ic I_{Q_0}}f(x):=\big( \sum_{\substack{ I \in \ic I_{Q_0} \\ I \subseteq Q_0}} \frac{|\langle f, \phi_I^1 \rangle|^2}{|I|} \one_I(x) \big)^\frac{1}{2}
\end{equation}
is constant on each $Q \in ch_{\ic S}(Q_0)$ and moreover, for every $x \in E_{Q_0}$,
\[
S_{\ic I_{Q_0}}f(x) \lesssim \frac{1}{|Q_0|^{\frac{1}{p_1}}} \Big\| \big( \sum_{\substack{ I \in \ic I \\ I \subseteq Q_0}} \frac{|\langle f, \phi_I^1 \rangle|^2}{|I|} \one_I \big)^\frac{1}{2}  \Big\|_{p_1}.
\]
The same inequality remains true outside of $E_{Q_0}$, by the definition \eqref{def-excp-set-sparse}. So that we have, for every $J_1 \in \ic I_{Q_0}$
\[
\frac{1}{|J_1|^\frac{1}{p_1}}   \Big\| \big(  \sum_{\substack{I \in \ic{I}_{Q_0} \\ I \subseteq J_1}} \frac{|\langle f, \phi_I^1  \rangle|^2}{|I|} \cdot \one_I  \big)^\frac{1}{2}   \Big\|_{p_1} \lesssim  \frac{1}{|Q_0|^\frac{1}{p_1}}\Big\| \big(  \sum_{\substack{I \in \ic{I} \\ I \subseteq Q_0}} \frac{|\langle f, \phi_I^1  \rangle|^2}{|I|} \cdot \one_I  \big)^\frac{1}{2}   \Big\|_{p_1}.
\]

Also, all $J_2 \in \ic I_{Q_0}$ intersect $\ds  \lbrace x \in Q_0: M( w \cdot \ci_{I_0})(x) > C \, \frac{1}{|Q_0|} \int_{\rr R^d} w \cdot \ci_{Q_0}(y) dy \rbrace ^c$, which implies 
\[
 \sup_{J_2 \in \ic{I}^+_{Q_0}} \frac{1}{|J_2|} \int_{\rr R^d} w \cdot \ci_{J_2} dx  \lesssim  \frac{1}{|Q_0|} \int_{\rr R^d} w \cdot \ci_{Q_0}  \, dx
\]

Using the subadditivity of $\| \cdot  \|_p^p$ and the result in Lemma \ref{lemma:localization-lemma}, we have
\begin{align*}
&\| \big( \sum_{I \in \mathcal{I}} \langle f, \phi_I^1  \rangle \,  \phi_I^2\big) \cdot w^\frac{1}{p} \|_p^p \lesssim \sum_{Q \in S} \| \big( \sum_{I \in \mathcal{I}_{Q}} \langle f, \phi_I^1  \rangle \,  \phi_I^2\big) \cdot w^\frac{1}{p} \|_p^p \\
& \lesssim \sum_{Q \in \ic S} \big( \sup_{J_1 \in \ic{I}_{Q}} \frac{1}{|J_1|^\frac{1}{p_1}}   \Big\| \big(  \sum_{\substack{I \in \ic{I}_{Q} \\ I \subseteq J_1}} \frac{|\langle f, \phi_I^1  \rangle|^2}{|I|} \cdot \one_I  \big)^\frac{1}{2}   \Big\|_{p_1}  \big)^p \, \big( \sup_{J_2 \in \ic{I}^+_{Q}} \frac{1}{|J_2|} \int_{\rr R^d} w \cdot \ci_{J_2} dx  \big) \cdot |Q| \\
&\lesssim \sum_{Q \in S} \big(  \frac{1}{|Q|^\frac{1}{p_1}}\Big\| \big(  \sum_{\substack{I \in \ic{I} \\ I \subseteq Q}} \frac{|\langle f, \phi_I^1  \rangle|^2}{|I|} \cdot \one_I  \big)^\frac{1}{2}   \Big\|_{p_1} \big)^p  \cdot \big(  \frac{1}{|Q|} \int_{\rr R^d} w \cdot \ci_Q  \, dx \big) \cdot |Q|.
\end{align*}

If $p>1$, we invoke a procedure that has already appeared in Proposition 20 of our previous \cite{BM3}. In this situation, we can use duality:
\[
\| \big( \sum_{I \in \mathcal{I}} \langle f, \phi_I^1  \rangle \,  \phi_I^2\big) \cdot w^\frac{1}{p} \|_p= \| \big( \sum_{I \in \mathcal{I}} \langle f, \phi_I^1  \rangle \,  \phi_I^2\big) \cdot w^\frac{1}{p} \, u\|_1,
\]
for some function $u \in L^{p'}$ with $\|u\|_{L^{p'}}=1$. Now we can apply the result of Theorem \ref{thm-sparse-dom} for $p=1$ to deduce the existence of a sparse collection $\ic S$ so that 
\[
\| \big( \sum_{I \in \mathcal{I}} \langle f, \phi_I^1  \rangle \,  \phi_I^2\big) \cdot w^\frac{1}{p} \, u \|_1 \lesssim \sum_{Q \in S} \big(  \frac{1}{|Q|^\frac{1}{p_1}}\Big\| \big(  \sum_{\substack{I \in \ic{I} \\ I \subseteq Q}} \frac{|\langle f, \phi_I^1  \rangle|^2}{|I|} \cdot \one_I  \big)^\frac{1}{2}   \Big\|_{p_1} \big)  \, \big( \frac{1}{|Q|} \int_{\rr R^d} w^\frac{1}{p} \, u \cdot \ci_Q  \, dx \big) \cdot |Q|.
\]

H\"older's inequality, first with respect to the measure $\ci_{Q} \, dx$ and with exponents $p+\epsilon$ and $(p+\epsilon)'$ yields
\begin{align*}
\| \big( \sum_{I \in \mathcal{I}} \langle f, \phi_I^1  \rangle \,  \phi_I^2\big) \cdot w^\frac{1}{p} \, u \|_1 & \lesssim \sum_{Q \in \ic S S} \big(  \frac{1}{|Q|^\frac{1}{p_1}}\big\| \big(  \sum_{\substack{I \in \ic{I} \\ I \subseteq Q}} \frac{|\langle f, \phi_I^1  \rangle|^2}{|I|} \cdot \one_I  \big)^\frac{1}{2}   \big\|_{p_1} \big) \\
&  \, \big( \frac{1}{|Q|} \int_{\rr R^d} w^\frac{p+\epsilon}{p} \cdot \ci_Q  \, dx \big)^\frac{1}{p+\epsilon}  \, \big( \frac{1}{|Q|} \int_{\rr R^d} u^{(p+\epsilon)'} \cdot \ci_Q  \, dx \big)^{\frac{1}{(p+\epsilon)'}} \cdot |Q|.
\end{align*}
Then we use again H\"older's inequality with respect to the discrete measure $\ds \ell^p(\ic S)$ to estimate the above expression by
\begin{align*}
& \Big( \sum_{Q \in S} \big(  \frac{1}{|Q|^\frac{1}{p_1}}\big\| \big(  \sum_{\substack{I \in \ic{I} \\ I \subseteq Q}} \frac{|\langle f, \phi_I^1  \rangle|^2}{|I|} \cdot \one_I  \big)^\frac{1}{2}   \big\|_{p_1} \big)^p \, \big( \frac{1}{|Q|} \int_{\rr R^d} w^\frac{p+\epsilon}{p} \cdot \ci_Q  \, dx \big)^\frac{p}{p+\epsilon} \cdot |Q|   \Big)^\frac{1}{p} \\
& \cdot \,  \Big( \sum_{Q \in S}  \big( \frac{1}{|Q|} \int_{\rr R^d} u^{(p+\epsilon)'} \cdot \ci_Q  \, dx \big)^{\frac{p'}{(p+\epsilon)'}} \cdot |Q| \Big)^\frac{1}{p'}.
\end{align*}

For the last term, use take advantage of the sparseness property, more exactly, we use the disjointness of the sets $\lbrace E(Q)  \rbrace_{Q \in \ic S}$:
\begin{align*}
&\Big( \sum_{Q \in S}  \big( \frac{1}{|Q|} \int_{\rr R^d} u^{(p+\epsilon)'} \cdot \ci_Q  \, dx \big)^{\frac{p'}{(p+\epsilon)'}} \cdot |Q| \Big)^\frac{1}{p'} \lesssim \Big( \sum_{Q \in S}  \big( \frac{1}{|Q|} \int_{\rr R^d} u^{(p+\epsilon)'} \cdot \ci_Q  \, dx \big)^{\frac{p'}{(p+\epsilon)'}} \cdot |E(Q)| \Big)^\frac{1}{p'} \\
&\lesssim \big\|  M_{(p+\epsilon)'} u   \big\|_{p'} \lesssim  \| u \|_{p'}=1.
\end{align*}
We are losing an $\epsilon$ (as small as we wish) in making sure that the maximal operator $M_{(p+\epsilon)'}$ is bounded on $L^{p'}$. We can choose $\epsilon$ so that $\ds \epsilon_p=\frac{p+\epsilon}{p}$. 
\end{proof}
\end{theorem}

Such a sparse estimate allows us to recover the weighted estimates from \cite{DingHanLu-weightedHardySpaces}, in the one-parameter case.

\begin{proposition}
\label{prop:weighted-est-SF}
Let $0<p<\infty$, $w \in A_\infty$ and $f$ a Schwartz function on $\rr R^d$; then
\begin{equation}
\label{eq:weighted-est}
\|f\|_{L^p(w)} \lesssim \|Sf\|_{L^p(w)}.
\end{equation} 
\begin{proof}
The weighted estimate follows easily once we prove a strengthening of the sparse estimate \eqref{eq:sparse-SF}: there exists a sparse collection of dyadic cubes $\ic S$ so that
\begin{equation}
\label{eq:sparse-w}
\| \big( \sum_{I \in \mathcal{I}} \langle f, \phi_I^1  \rangle \,  \phi_I^2\big) \cdot w^\frac{1}{p} \|_p^p \lesssim \sum_{Q \in S} \big(  \frac{1}{|Q|^\frac{1}{p_1}}\Big\| \big(  \sum_{\substack{I \in \ic{I} \\ I \subseteq Q}} \frac{|\langle f, \phi_I^1  \rangle|^2}{|I|} \cdot \one_I  \big)^\frac{1}{2}   \Big\|_{p_1} \big)^p w(E(Q)).
\end{equation}

If such an estimate were true, we could deduce that
\begin{align*}
\| \big( \sum_{I \in \mathcal{I}} \langle f, \phi_I^1  \rangle \,  \phi_I^2\big) \|_{L^p(w)}^p \lesssim \sum_{Q \in S}  \big( \inf_{y \in Q} M ( |S \,f |^{p_1} )(y)  \big)^\frac{p}{p_1}  \cdot w(E(Q)),
\end{align*}
and in consequence, 
\[
\|  \sum_{I \in \mathcal{I}} \langle f, \phi_I^1  \rangle \,  \phi_I^2 \|_{L^p(w)}^p \lesssim \int_{\rr R^d} M ( |S \,f |^{p_1} )(x)  \big)^\frac{p}{p_1} \, w(x) \, dx.
\]

So far, no information was required on $p_1$; it suffices to choose $p_1<p$ and so that $w \in A_{\frac{p}{p_1}}$ (this will assure that $ M$ is bounded on $L^{\frac{p}{p_1}}(w)$) to obtain that 
\[
\| \sum_{I \in \mathcal{I}} \langle f, \phi_I^1  \rangle \,  \phi_I^2  \|_{L^p(w)} \lesssim \|Sf\|_{L^p(w)}.
\]
This is possible since $\ds w \in A_\infty=\bigcup_{q>1} A_q$. The final inequality \eqref{eq:weighted-est} is deduced thanks to formula \eqref{desc}.

We are left with showing how \eqref{eq:sparse-SF} implies \eqref{eq:sparse-w}. We recall that 
{\fontsize{10}{10}\[
\big( \frac{1}{|Q|} \int_{\rr R^d} w^{1+\epsilon_p} \, \ci_Q  \, dx \big)^\frac{1}{1+\epsilon_p} \leq \sum_{\ell \geq 0} 2^{- \ell M} \big( \frac{1}{|Q|} \int_{2^\ell Q} w^{1+\epsilon_p}   \, dx \big)^\frac{1}{1+\epsilon_p} \leq \sum_{\ell \geq 0} 2^{- \ell M} 2^\frac{\ell d}{1+\epsilon_p} \big( \frac{1}{|2^\ell Q|} \int_{2^\ell Q} w^{1+\epsilon_p}   \, dx \big)^\frac{1}{1+\epsilon_p}.
\]}

Now we use the \emph{Reverse H\"older} property of the weight $w$: there exists $\epsilon_w$ so that 
\[
\big( \frac{1}{|2^\ell Q|} \int_{2^\ell Q} w^{1+\epsilon_w}   \, dx \big)^\frac{1}{1+\epsilon_w} \lesssim \frac{1}{|2^\ell Q|} \int_{2^\ell Q} w \, dx.
\]

If we pick $\epsilon_p < \epsilon_w$, then the $L^{1+\epsilon_p}$ average in \eqref{eq:sparse-SF} can be replaced by an $L^1$ average (note that, for $0<p \leq 1$, we have from the start $\epsilon_p=0$). Hence, we have 
{\fontsize{10}{10}\begin{align*}
\| \big( \sum_{I \in \mathcal{I}} \langle f, \phi_I^1  \rangle \,  \phi_I^2\big) \cdot w^\frac{1}{p} \|_p^p \lesssim \sum_{\ell \geq 0} 2^{- \ell M} 2^\frac{\ell d}{1+\epsilon_p} \sum_{Q \in S} \big(  \frac{1}{|Q|^\frac{1}{p_1}}\Big\| \big(  \sum_{\substack{I \in \ic{I} \\ I \subseteq Q}} \frac{|\langle f, \phi_I^1  \rangle|^2}{|I|} \cdot \one_I  \big)^\frac{1}{2}   \Big\|_{p_1} \big)^p  2^{-\ell d} w(2^\ell Q).
\end{align*}}

All we need to do is compare $w(2^\ell Q)$ and $w(E(Q))$. We know that $|Q|<2 \,|E(Q)|$ and $w \in A_{\infty}$. Then $w \in A_{q_w}$ for some $q_w >1$ and in consequence (see inequality (7.2) of \cite{Duoa-book})
\[ 
w(2^\ell Q) \big( \frac{|E(Q)|}{|2^\ell Q|} \big) \lesssim w(E(Q)) \quad \Longleftrightarrow \quad w(2^\ell Q) \lesssim 2^{\ell \, d\, q_w} w(E(Q)).
\]

If $M$, the decaying exponent of the auxiliary weights $\ci_{Q}$ (see Definition \ref{def:ci}) satisfies $\ds d \, q_w < M$, then we can sum in $\ell \geq 0$ and we are done. Since $M$ can be as large as we wish, we can arrange for this condition to be satisfied.
\end{proof}
\end{proposition}

We note that the sparse domination result \eqref{eq:sparse-SF} of Theorem \ref{thm-sparse-dom} implies, for any collection $\ic I$ of dyadic squares and any fixed dyadic square $I_0$:
\begin{equation}
\label{eq:sparse-max-weight-only}
\| \big( \sum_{I \in \mathcal{I}(I_0)} \langle f, \phi_I^1  \rangle \,  \phi_I^2\big) \|_{L^p(w)}^p \lesssim \big( \sup_{J_2 \in \ic I^+(I_0)} \frac{1}{|J_2|} \int_{\rr R^d} w(x) \, \ci_{J_2} dx  \big) \, \big\| \ic S_{\ic I(I_0)} f  \big\|_p^p.
\end{equation}

This observation will be useful shortly, as we will show that it is possible to prove a multiple vector-valued, weighted result without making use of extrapolation. 
\begin{proposition}
\label{prop-mvv-weighted-est}
Let $0<p<\infty, \, 0<Q<\infty$ and $w \in A_{\infty}$; then for any $L^Q$-valued Schwartz function $f$ on $\rr R^d$, we have
\[
\|f\|_{L^p(L^Q; d\, w)} \lesssim \|S \,f\|_{L^p(L^Q; d\, w)}.
\]
\end{proposition} 

The proof combines together all the previous techniques used for deducing multiple vector-valued estimates in Section \ref{sec:proof-discretized-main-thm} and weighted estimates. We sketch the proof of the crucial maximal inequality (the equivalent of \eqref{eq:local-eq} of Lemma \ref{lemma:localization-lemma}) in the case of $\ell^q$-valued functions, where $q<1$. The case $q \geq 1$ is in fact easier, since duality is available. The general multiple vector-valued case, corresponding to a general $n$-tuple $Q$, follows by induction over $n$.

\begin{lemma}
\label{lemma:local-est-vv}
Let $0<q<1$ and $0<p\leq q$; let $\ic I$ be a finite collection of dyadic squares in $\rr R^d$, $I_0$ a fixed dyadic square, $f : \rr R^d \to \rr C$ a Schwartz function and $w$ a locally integrable, positive function. Then for any $0<p_1<\infty$,
{\fontsize{10}{10}\begin{align}
\label{eq:ineq-max-vv-w}
&\big\| \big( \sum_k \big| \sum_{I \in \mathcal{I}(I_0)} \langle f_k, \phi_I^1 \rangle \, \phi_I^2 \big|^q \big)^\frac{1}{q} \big\|_{L^p(w)} \\
& \lesssim \big( \sup_{J_1 \in \ic{I}(I_0)} \frac{1}{|J_1|^\frac{1}{p_1}}   \Big\| \Big( \sum_k | \big(  \sum_{\substack{I \in \ic{I}(I_0) \\ I \subseteq J_1}} \frac{|\langle f_k, \phi_I^1  \rangle|^2}{|I|} \cdot \one_I  \big)^\frac{q}{2} \Big)^\frac{1}{q}   \Big\|_{p_1}  \big) \, \big( \sup_{J_2 \in \ic{I}^+(I_0)} \frac{1}{|J_2|} \int_{\rr R} w \cdot \ci_{J_2} dx  \big)^\frac{1}{p} \cdot |I_0|^\frac{1}{p}, \nonumber
\end{align}}with an implicit constant independent of the collection $\ic I$ and of the functions $f$ and $w$.
\begin{proof}

We note that $\| \cdot \|^p_{L^p(\ell^q; d\, w)}$ is subadditive, and hence, using the decomposition \eqref{desco},
\[
\big\| \big( \sum_k \big| \sum_{I \in \mathcal{I}(I_0)} \langle f_k, \phi_I^1 \rangle \, \phi_I^2 \big|^q \big)^\frac{1}{q} \big\|_{L^p(w)}^p \lesssim \sum_{\ell \geq 0} 2^{-\ell \, p\, M} \big\| \big( \sum_k \big| \sum_{I \in \mathcal{I}(I_0)} \langle f_k, \phi_I^1 \rangle \, \tilde{\phi}_{I, \ell}^2 \big|^q \big)^\frac{1}{q} \big\|_{L^p(w)}^p.
\]

Since $p \leq q$ and all the functions $\tilde \phi_{I, \ell}^2$ are supported inside $2^\ell I_0$:
\[
\big\| \big( \sum_k \big| \sum_{I \in \mathcal{I}(I_0)} \langle f_k, \phi_I^1 \rangle \, \tilde{\phi}_{I, \ell}^2 \big|^q \big)^\frac{1}{q} \big\|_{L^p(w)} \lesssim \big\| \big( \sum_k \big| \sum_{I \in \mathcal{I}(I_0)} \langle f_k, \phi_I^1 \rangle \, \tilde{\phi}_{I, \ell}^2 \big|^q \big)^\frac{1}{q} \big\|_{L^q(w)} \cdot \|  \one_{2^\ell \, I_0}  \|_{L^{\tau}(w)},
\]
where $\frac{1}{p}=\frac{1}{q}+\frac{1}{\tau}$. For the first term on the right hand side, we use Fubini and the known scalar version of Lemma \ref{lemma:local-est-vv} (more precisely, inequality \eqref{eq:sparse-max-weight-only} above):
{\fontsize{9}{9}\begin{align*}
& \big\| \big( \sum_k \big| \sum_{I \in \mathcal{I}(I_0)} \langle f_k, \phi_I^1 \rangle \, \tilde{\phi}_{I, \ell}^2 \big|^q \big)^\frac{1}{q} \big\|_{L^q(w)}^q =\sum_k  \big\|  \sum_{I \in \mathcal{I}(I_0)} \langle f_k, \phi_I^1 \rangle \, \tilde{\phi}_{I, \ell}^2 \big\|_{L^q(w)}^q \lesssim \sum_k  \big\| \ic S_{\ic I(I_0)} f_k  \big\|_q ^q \, \big( \sup_{J_2 \in \ic I^+(I_0)} \frac{1}{|J_2|} \int_{\rr R^d} w(x) \, \ci_{J_2} dx  \big) \\
& \lesssim \Big( \frac{1}{|I_0|^\frac{1}{q}}  \big\| \big( \sum_k  |\ic S_{\ic I(I_0)} f_k|^q \big)^\frac{1}{q} \big\|_q   \Big)^q \cdot \big( \sup_{J_2 \in \ic I^+(I_0)} \frac{1}{|J_2|} \int_{\rr R^d} w(x) \, \ci_{J_2} dx  \big) \cdot |I_0| 
\end{align*}}

By a vector-valued version of John-Nirenberg's inequality, which was also used in proving the multiple vector-valued version of Theorem \ref{mainth}, the above can be estimated by
\[
\Big( \sup_{J_1 \in \ic I(I_0)}  \frac{1}{|I_0|^\frac{1}{p_1}}  \big\| \big( \sum_k  |\ic S_{\ic I(I_0)} f_k|^q \big)^\frac{1}{q} \big\|_{p_1}   \Big)^q \cdot \big( \sup_{J_2 \in \ic I^+(I_0)} \frac{1}{|J_2|} \int_{\rr R^d} w(x) \, \ci_{J_2} dx  \big) \cdot |I_0|, 
\]
where $0<p_1<\infty$ is any Lebesgue exponent.

On the other hand, 
\[ 2^{-\ell p M/2}  \|  \one_{2^\ell \, I_0}  \|_{L^{\tau}(w)}^p \lesssim \big(  \frac{1}{|I_0|} \int_{\rr R^d} w(x) \, \ci_{I_0} dx  \big)^\frac{p}{\tau} |I_0|^\frac{p}{\tau}.
\]
After summing in $\ell \geq 0$, we get the inequality \eqref{eq:ineq-max-vv-w}.

\end{proof}
\end{lemma}

Applying the usual stopping time, the maximal inequality of Lemma \ref{lemma:local-est-vv} will imply a vector-valued version of Theorem \ref{thm-sparse-dom}. We leave the details to the interested reader. Although Lemma \ref{lemma:local-est-vv} is stated for $p \leq q$, a vector-valued version of Theorem \ref{thm-sparse-dom} is valid for any Lebesgue exponents, as we can pass from lower Lebesgue exponents to larger ones at the expense of loosing an $\epsilon$.

\subsubsection{The multi-parameter case}
\label{sec:weighted-ineq-multi-param}
The multi-parameter version of Proposition \ref{prop:weighted-est-SF} follows easily from the properties of the weights $A_{\infty, Rectangle}(\rr R^{d_1} \times \ldots \times \rr R^{d_N})$. We will only illustrate the scalar bi-parameter case, but state the result in its generality.

\begin{proposition}
\label{prop:mvv-multi-param-weighted}
Let $0<p<\infty$, $0<Q<\infty$; then for any $w \in A_{\infty, Rectangle}(\rr R^{d_1} \times \ldots \times \rr R^{d_N})$ and any $L^Q$-valued Schwartz function $f$, 
\[
\| f\|_{L^p(L^Q)(w)} \leq C \,\|S_{d_1} \otimes\ldots \otimes  S_{d_N} ( f)\|_{L^p(L^Q)(w)}.
\]
\begin{proof}
In fact, we will prove that 
\begin{equation}
\label{eq:bi-param-weight}
\| f\|_{L^p(w)} \leq C\|S_{d_1} \otimes  S_{d_2} ( f)\|_{L^p(w)},
\end{equation}
for any $w \in A_{\infty, Rectangle}(\rr R^{d_1} \times \rr R^{d_2})$. An important property of the weights in the class $A_{\infty, Rectangle}(\rr R^{d_1} \times \ldots \times \rr R^{d_N})$ is that if we fix one of the variables, we still obtain an $A_\infty$ weight in the other variable and we can use the one-parameter result:
{\fontsize{10}{10}\begin{equation}
\label{eq:weights-fix-var}
w_y(x)=w(x, y) \in A_{\infty}(\rr R^{d_1}) \text{  for a.e.  } y \in \rr R^{d_2}, \quad w_x(y)=w(x, y) \in A_{\infty}(\rr R^{d_2}) \text{   for a.e.  } x \in \rr R^{d_2}.
\end{equation}}

We start by fixing the variable $y$ ; then $f_y(x):=f(x,y)$ is a function on $\rr R^{d_1}$. By Proposition \ref{prop:weighted-est-SF},
\[
\int_{\rr R^{d_1}} |f(x, y)|^p w(x, y) dx \lesssim \int_{\rr R^{d_1}} |S_{d_1} f_y (x)| dx  = \int_{\rr R^{d_1}}  \big( \sum_{k} |Q_k(f_y)(x)|^2   \big)^\frac{p}{2} w(x, y) dx.
\]
Above,
\begin{equation}
\label{eq:about-notation-Q_k}
Q_k(f_y)(x):=Q_k^1 f(x, y):= Q_k^x f(y):=\int_{\rr R^{d_1}} f(x-s, y) \psi_k(s) ds.
\end{equation}

If we integrate with respect to $y$ and use Fubini, we have
\[
\int_{\rr R^{d_1}} \int_{\rr R^{d_2}}  |f(x, y)|^p w(x, y) dy dx \lesssim  \int_{\rr R^{d_1}} \int_{\rr R^{d_2}}   \big( \sum_{k} |(Q_k^x f)(y)|^2   \big)^\frac{p}{2} w(x, y)dy\, dx.
\]
Now we consider $x$ fixed and we apply Proposition \ref{prop:weighted-est-SF} (or more specifically an $\ell^2$-valued extension which follows also from a well-known result of Marcinkiewicz and Zygmund \cite{MarcinkiewiczZygmundL2Extensions}) to the sequence of functions $(Q_k^x f)_{k \in \rr Z}$:
{\fontsize{9}{9}\[
\int_{\rr R^{d_2}}   \big( \sum_{k} |(Q_k^x f)(y)|^2   \big)^\frac{p}{2} w(x, y)\,dy \lesssim \int_{\rr R^{d_2}}  \big( \sum_{k} |S_{d_2}(Q_k^x f)(y)|^2   \big)^\frac{p}{2} w(x, y)\,dy= \int_{\rr R^{d_2}}  \big( \sum_{k} |\sum_l Q_l(Q_k^x f)(y)|^2   \big)^\frac{p}{2} w(x, y)\,dy.
\]}
Here it is useful that we can interchange the role played by the variables: if $x$ is fixed, $w(x, \cdot)$ is still an $A_\infty$ weight and vice-versa.

We need to understand the last expression:
{\fontsize{10}{10}
\[
 Q_l(Q_k^x f)(y)= \int_{\rr R^{d_2}} (Q_k^x f)(y-t) \, \psi_l(t) dt = \int_{\rr R^{d_2}} \big(\int_{\rr R^{d_1}} f(x-s, y-t) \psi_k(s) ds\big)  \psi_l(t) dt= f \ast (\psi_k \otimes \psi_l)(x, y),
\]}
so that 
\[
\big( \sum_{k} |S_{d_2}(Q_k^x f)(y)|^2   \big)^\frac{p}{2} = \big( \sum_{k} |\sum_l Q_l(Q_k^x f)(y)|^2   \big)^\frac{1}{2} =S_{d_1} \otimes S_{d_2}(f)(x,y).
\]
Integrating in $x$ we obtain \eqref{eq:bi-param-weight}. Note that here it is important that we can use Fubini, fix one of the variable and perform the usual one-parameter analysis; in particular, the properties \eqref{eq:weights-fix-var} are critical. For mixed-norm estimates most of the weighted results are known only for weights that tensorize: $w(x,y)=u(x)\, v(y)$, the reason being that Fubini and property \eqref{eq:weights-fix-var} do not hold any longer.
\end{proof}
\end{proposition}

\end{document}